\documentclass[3p]{elsarticle}

\usepackage{makeidx}
\usepackage{amsfonts}
\usepackage{amsmath}
\usepackage{amssymb}
\usepackage{amscd}
\usepackage{xspace}
\usepackage{color}
\usepackage{array,amssymb,amsmath,amsbsy}
\usepackage{enumerate}
\usepackage[bottom]{footmisc}
\usepackage[all,cmtip]{xy}
\usepackage{hyperref}
\usepackage{wrapfig}
\usepackage{simplemargins}

\setallmargins{1in}
\usepackage[colorinlistoftodos,backgroundcolor=blue!10,linecolor=red,bordercolor=red]{todonotes}

\newtheorem{theorem}{Theorem}[section]
\newtheorem{lemma}[theorem]{Lemma}

\newtheorem{proposition}[theorem]{Proposition}
\newtheorem{corollary}[theorem]{Corollary}

\newproof{proof}{Proof}
\newproof{claim_proof}{Proof of Claim}

\def\Aut{{\rm Aut}}
\def\Autbo#1{{\rm Aut}_{\bo #1}(#1)}
\def\Stab{{\rm Stab}}

\def\Fix{{\rm Fix}}
\def\Ker{{\rm Ker}}
\def\Iso{{\rm Iso}}
\def\O{\calO}
\def\id{{\rm id}}
\def\bo{\partial} % boundary
\def\int{\mathring} % interior

\def\planar{\sf{PLANAR}}
\def\connectedplanar{\hbox{\rm \sffamily connected \planar}\xspace}
\def\R{\mathbb R}

\def\bV{\boldsymbol{V}}
\def\bE{\boldsymbol{E}}
\def\bD{\boldsymbol{D}}
\def\bv{\boldsymbol{v}}
\def\be{\boldsymbol{e}}
\def\bd{\boldsymbol{d}}

\def\rcover{\textsc{RegularCover}}
\def\cover{\textsc{$H$-Cover}}

\newenvironment{packed_enum}{
	\begin{enumerate}
		\setlength{\itemsep}{1pt}
	    \setlength{\parskip}{0pt}
		\setlength{\parsep}{0pt}
}{\end{enumerate}}

\newenvironment{packed_itemize}{
	\begin{itemize}
		\setlength{\itemsep}{1pt}
	    \setlength{\parskip}{0pt}
		\setlength{\parsep}{0pt}
}{\end{itemize}}

\newenvironment{packed_head_enum}[1]{
	\begin{enumerate}[#1]
		\setlength{\itemsep}{1pt}
	    \setlength{\parskip}{0pt}
		\setlength{\parsep}{0pt}
}{\end{enumerate}}

\newcommand{\heading}[1]{\medskip\par\noindent{\bf #1}}

\def\computationproblem#1#2#3#4{% {problem_name}{input}{output}
	\vskip 1ex
	\begin{center}
	\fbox{\begin{tabular}{rp{#4}}
	{\bf Problem:\enspace}&#1\\
	{\bf Input:\enspace}&#2\\
	{\bf Question:\enspace}&#3\\
	\end{tabular}}
	\end{center}
	\vskip 1ex
}

\def\gS{\mathbb{S}} \def\gC{\mathbb{C}} \def\gA{\mathbb{A}} \def\gD{\mathbb{D}}

\def\calA{{\cal A}}

  \def\calO{{\cal O}} \def\calP{{\cal P}}

\def\cNP{\hbox{\rm \sffamily NP}}

\begin{document}
\title{3-connected Reduction for Regular Graph Covers\tnoteref{conference,support}}
\tnotetext[conference]{This paper continues the research started in ICALP 2014~\cite{fkkn} and
extends its results. For a unified description of the results of this and the follow-up papers, see
the PhD thesis~\cite{phd_thesis}. For a structural diagram visualizing our results, see
\url{http://pavel.klavik.cz/orgpad/regular_covers.html} (supported for Firefox and Google Chrome).}
\tnotetext[support]{This work was initiated during workshops Algebraic, Topological and Complexity
Aspects of Graph Covers (ATCACG). The authors are supported by CE-ITI (P202/12/G061 of GA\v{C}R).
The first author is also supported by the project Kontakt LH12095, the second, the third, and the
fourth authors by Charles University as GAUK 1334217, the fourth author by the Ministry of Education
of the Slovak Republic, the grant VEGA 1/0487/17, and by the project APVV-15-0220 of Slovak Research
and Development Agency.}

\author[cunidam]{Ji\v{r}\'i Fiala}
\ead{fiala@kam.mff.cuni.cz}
\author[cunidam,cunicsi]{Pavel Klav\'ik}
\ead{klavik@iuuk.mff.cuni.cz}
\author[cunidam]{Jan Kratochv\'il}
\ead{honza@kam.mff.cuni.cz}
\author[mbu,ntis]{Roman Nedela}
\ead{nedela@savbb.sk}

\address[cunidam]{Department of Applied Mathematics, Faculty of Mathematics and
	   	Physics, Charles University, Prague, Czech Republic.}
\address[cunicsi]{Computer Science Institute, Faculty of Mathematics and Physics, Charles
		University, Prague, Czech Republic.}
\address[mbu]{Institute of Mathematics and Computer Science SAS, Bansk\'a Bystrica, Slovak republic.}
\address[ntis]{European Centre of Excellence NTIS, University of West Bohemia, Pilsen, Czech Republic.}

\begin{abstract}
A graph $G$ \emph{covers} a graph $H$ if there exists a locally bijective homomorphism from $G$ to
$H$. We deal with \emph{regular coverings} in which this homomorphism is prescribed by an action of a
semiregular subgroup $\Gamma$ of $\Aut(G)$; so $H \cong G / \Gamma$.  In this paper, we study
the behaviour of regular graph covering with respect to 1-cuts and 2-cuts in $G$.

We describe \emph{reductions} which produce a series of graphs $G = G_0,\dots,G_r$ such that
$G_{i+1}$ is created from $G_i$ by replacing certain inclusion minimal subgraphs with colored edges.
The process ends with a \emph{primitive graph} $G_r$ which is either 3-connected, or a cycle, or
$K_2$. This reduction can be viewed as a non-trivial modification of reductions of Mac Lane (1937),
Trachtenbrot (1958), Tutte (1966), Hopcroft and Tarjan (1973), Cuningham and Edmonds (1980), Walsh
(1982), and others. A novel feature of our approach is that in each step all essential information
about symmetries of $G$ are preserved.

A regular covering projection $G_0\to H_0$ induces regular covering projections $G_i \to H_i$ where
$H_i$ is the $i$-th quotient reduction of $H_0$. This property allows to construct all possible
quotients $H_0$ of $G_0$ from the possible quotients $H_r$ of $G_r$. By applying this method to
planar graphs, we give a proof of Negami's Theorem (1988). Our structural results are also used in
subsequent papers for regular covering testing when $G$ is a planar graph and for an inductive
characterization of the automorphism groups of planar graphs (see Babai (1973) as well).

\end{abstract}


\maketitle

\section{Introduction}

\begin{figure}[t!]
\centering
\includegraphics{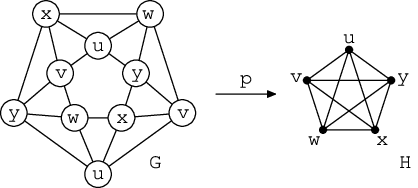}
\caption{A regular covering projection $p$ from a graph $G$ to one of its quotients $H$. For every
vertex $v \in \bV(G)$, the image $p(v)$ is written in the circle.}
\label{fig:big_picture}
\end{figure}

The notion of \emph{covering} originates in topology to describe local similarity of two topological
spaces.  In this paper, we study coverings of graphs in a more restricting version called
\emph{regular covering}, for which the covering projection is described by a semiregular action of a
group; see Section~\ref{sec:preliminaries} for the formal definition. If $G$ regularly covers $H$,
then $H$ is called a \emph{regular quotient} of $G$, or just a \emph{quotient}. See
Fig.~\ref{fig:big_picture} for an example.

Regular graph covers have many applications in graph theory, for instance they were used to solve
the Heawood map coloring problem~\cite{ringel_youngs,gross_tucker} and to construct arbitrarily large
highly symmetrical graphs~\cite{biggs}. The concept of a regular covering of graphs gives rise to a
powerful construction of large graphs with prescribed properties from smaller
ones~\cite{hoffman_construction,siagiova,lifts_maria}.  It can be demonstrated by the well-known
construction of a Cayley graph, where a large graph is defined by specifying few generators of a
group. While each Cayley graph can be viewed as regular cover over a one-vertex graph, a regular
cover is a generalization of the construction of a Cayley graph, where the one-vertex quotient is
replaced by an arbitrary connected graph.

\subsection{Our Results}

In this paper, we fully describe the behaviour of regular covering with respect to 1-cuts and 2-cuts
in $G$, for the missing definitions the reader is reffered to Section~\ref{sec:preliminaries}.
Since any regular covering is equivalent to the natural projection $G\to G/\Gamma$, where
$\Gamma\leq \Aut(G)$ is semiregular on the vertices and darts of $G$, but not necessarily on the
edges of $G$, the problem  is closely related to an investigation of the behaviour of 1- and 2-cuts
in a semiregular action of a subgroup of the automorphism group. Since every graph automorphism
fixes the central block of the block decomposition, to investigate the action of $\Gamma$ on the set
of articulations is not difficult. In particular, 1-cut is mapped onto a 1-cut in a (regular)
covering. However, the behaviour of a regular covering on 2-cuts is complex. Our main result,
Theorem~\ref{thm:quotient_expansion} describes all possible quotients of some graph class provided
we understand quotients of 3-connected graphs in this class. Our result applies to the class of
planar graphs, since the 3-connected planar graphs are 1-skeletons of polehedra, and the quotients
are given by semiregular actions of spherical groups on the polyhedra, see
Sections~\ref{sec:algo_cons} and~\ref{sec:structural_cons} for details.

Let us explain briefly our approach. We process the input graph $G$ by a series of \emph{reductions}
replacing some subgraphs of $G$, called \emph{atoms}, separated by 1- and 2-cuts, by edges.  This
natural idea of the reduction was first introduced in the seminal papers by Mac Lane~\cite{maclane}
and Trakhtenbrot~\cite{trakhtenbrot}. It was further extended
in~\cite{tutte_connectivity,quadratic_isomorphism_planar,hopcroft_tarjan_dividing,cunnigham_edmonds,walsh,bienstock}
and studied for infinite graphs in~\cite{droms1995structure}.  This decomposition can be represented
by a tree whose nodes are 3-connected graphs, and this tree is known in the literature mostly under
the name \emph{SPQR tree}~\cite{spqr1,spqr2,spqr3,spqr_linear}.  Our reduction has two key
differences:
\begin{packed_itemize}
\item The aforementioned papers apply the reduction exclusively to 2-connected graphs. In contrast,
in each step we simultaneously reduce subgraphs separated both by 1-cuts and 2-cuts. This requires
to introduce the definition of atoms in a proper way. Since a quotient of a 2-connected graph might
have 1-cuts, such a unified treatment is desirable.
\item The reduction is augemented by colored edges (encoding different isomorphism classes) of three
different types (encoding different symmetry types of atoms). This allows to capture the changes in
the automorphism group and its semiregular subgroups.
\item Since we allow non-trivial edge-stabilisers in the action of a semiregular group on a graph,
the respective edge-orbit is mapped by the covering projection onto an ``edge'' incident to exactly one
vertex and not being a loop. Such an ``edge'' is called a half-edge. Of course, allowing the
existence of half-edges requires to extend the usual definition of a graph. As a result we obtain
more general statements, than in the frame of the classical theory of graphs.
\end{packed_itemize}

\heading{Atoms and Reductions.} In Section~\ref{sec:atoms}, we introduce the essential definition of
an \emph{atom}. The atoms are, roughtly speaking, inclusion-minimal subgraphs with respect to 1-cuts
and 2-cuts which cannot be further simplified. They are essentially paths, cycles, stars, dipoles or
3-connected.  The reduction constructs a series of graphs $G = G_0,G_1,\dots,G_r$.  The reduction
from $G_i$ to $G_{i+1}$ is done by replacing all the atoms of $G_i$ by colored edges, where the
colors encode the isomorphism classes of atoms. The last (irreducible) graph in the sequence,
denoted by $G_r$, is called \emph{primitive}.  It is either very simple ($K_1$, $K_2$ or a cycle),
or it is 3-connected; or $G_r$ is obtained from these graphs by attaching single pendant edges to
some of the vertices.  Following the literature we call the reduction process the
\emph{$3$-connected reduction}.

When the graph $G$ is not 3-connected, we consider its block-tree. The central block plays the key
role in every regular covering projection. The reason is that a covering $G\to H$ behaves
non-trivially only on the central block; the remaining blocks are mapped by the covering onto the
isomorphic copies in $H$. Therefore the atoms are defined with respect to the central block.  We
distinguish three types of atoms:
\begin{packed_itemize}
\item \emph{Proper atoms} are inclusion-minimal subgraphs separated by a 2-cut inside a
block.
\item \emph{Dipoles} are formed by the sets of all parallel edges joining two vertices.
\item \emph{Block atoms} are blocks which are leaves of the block-tree, or stars consisting of all
pendant edges attached to a vertex. The central block is never a block atom.
\end{packed_itemize}

The reduction from $G_i$ to $G_{i+1}$ is defined in a way that there exists an induced
\emph{reduction epimorphism} $\Phi_i : \Aut(G_i)\to \Aut(G_{i+1})$ which possesses some nice
properties; see Proposition~\ref{prop:reduction_homomorphism} for details. Using it, we can describe
the respective change of the automorphism group explicitly:

\begin{proposition} \label{prop:group_extension}
If $G_i$ is reduced to $G_{i+1}$, then
$$\Aut(G_{i+1}) \cong \Aut(G_i) / \Ker(\Phi_i).$$
\end{proposition}

\heading{Expansions.} We aim to investigate how the knowledge of regular quotients of $G_{i+1}$ can
be used to construct all regular quotients of $G_i$. To do so, we introduce the reversal of the
reduction called the \emph{quotient expansion}. If $H_{i+1} = G_{i+1} / \Gamma_{i+1}$, then the
quotient expansion produces $H_i$ by replacing colored edges back by atoms. To do this, we have to
understand how regular covering behaves with respect to the atoms. Inspired by Negami~\cite{negami},
we show that each proper atom/dipole has three possible types of quotients that we call an
\emph{edge-quotient}, a \emph{loop-quotient} and a \emph{half-quotient}. The edge-quotient and the
loop-quotient are uniquely determined but an atom may have many non-isomorphic half-quotients.

\begin{figure}[b!]
\centering
\includegraphics{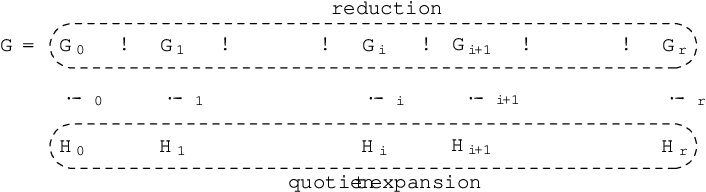}
\caption{The reduction is on top, the quotient expansion is on bottom. It holds that $H_i = G_i /
\Gamma_i$ and $\Gamma_i$ is a group extension of $\Gamma_{i+1}$.}
\label{fig:reduction_expansion_diagram}
\end{figure}

The constructed quotients contain colored edges, loops and half-edges corresponding to atoms.  Each
half-edge in $H_{i+1}$ is an image of a half-edge, or of a halvable edge $e$ in $G_{i+1}$ such that
an automorphism of $\Gamma_{i+1}$ fixes $e$ and swaps the vertices incident with $e$, see the next
section for exact definitions.  The following theorem is our main result, it describes every
possible expansion of $H_{i+1}$ to $H_i$:

\begin{theorem} \label{thm:quotient_expansion}
Let $G = G_0,\dots,G_r$ be the reduction series for a graph $G$.  Then every quotient $H_i$ of
$G_i$, for $i \in \{0,\dots,r-1\}$, can be constructed from some quotient $H_{i+1}$ of $G_{i+1}$ by
replacing each edge, loop and half-edge of $H_{i+1}$ by the subgraph corresponding to the edge-, the
loop-, or a half-quotient  of an atom of $G_i$, respectively. 
\end{theorem}

Suppose that some regular quotient of the primitive graph $G_r$ is chosen, so $H_r = G_r /
\Gamma_r$. The above theorem allows to describe all regular quotients $H$ of $G$ rising from $H_r$,
as depicted in the diagram in Fig.~\ref{fig:reduction_expansion_diagram}.

\subsection{Algorithmic and Complexity Consequences} \label{sec:algo_cons}

Our main algorithmic motivation is the study of the computational complexity of regular covering
testing:
\computationproblem
{\rcover}
{Connected graphs $G$ and $H$.}
{Does $G$ regularly cover $H$?}
{4.5cm}

\noindent Our structural results have the following algorithmic implications, described
in~\cite{fkkn,fkkn16}:

\begin{theorem}[Fiala et al.~\cite{fkkn,fkkn16}] \label{thm:algorithmic}
If $G$ is planar, we can solve \rcover\ in time $\O(\bv(G)^c \cdot 2^{\be(H)/2})$, where $c$ is
a constant, $\bv(G)$ is the number of vertices of $G$, and $\be(H)$ is the number of edges of $H$.
\end{theorem}

Theorem~\ref{thm:quotient_expansion} suggests that there might be exponentially many quotients of
$G$, and so this algorithm has to test efficiently whether $H$ is one of them. However, for
every fixed graph $H$, the constructed algorithm runs in polynomial time.

\heading{Relations to General Covers.} The aforementioned decision problem is closely related to the
complexity of general covering testing which was widely studied before. We try to understand how
much the additional algebraic structure influences the computational complexity. Study of the
complexity of general covers was pioneered by Bodlaender~\cite{bodlaender} in the context of
networks of processors in parallel computing, and for fixed target graph was first asked by Abello
et al.~\cite{AFS}. The problem \cover\ asks whether an input graph $G$ covers a fixed graph $H$. The
general complexity is still unresolved but papers~\cite{kratochvil97,fiala00} show that it is
\cNP-complete for every $r$-regular graph $H$ where $r \ge 3$. For a survey, see~\cite{FK}.

The complexity results concerning graph covers are mostly \cNP-complete. In our impression, the
additional algebraic structure of regular graph covers makes the problem easier, as shown by the
following two contrasting results. The problem \cover\ remains \cNP-complete for several small
fixed graphs $H$ (such as $K_4$, $K_5$) even for planar inputs $G$~\cite{planar_covers}. On the
other hand, Theorem~\ref{thm:algorithmic} shows that for planar graphs $G$ the problem \rcover\ is
fixed-parameter tractable in the number of edges of $H$.

\heading{Relations to Cayley Graphs and Graph Isomorphism.}
The notion of regular graph covers builds a bridge between two seemingly different problems. If the
base graph $H$ is a one-vertex graph, it corresponds to the problem of recognition of Cayley graphs
whose complexity is not known. A polynomial-time algorithm is known only for circulant
graphs~\cite{circulant_recognition}. When both graphs $G$ and $H$ have the same size, we get graph
isomorphism testing. Our results are far from solving these problems, but we believe that better
understanding of \rcover\ can also shed new light on these famous problems.

Theoretical motivation for studying graph isomorphism is very similar to \rcover. For practical
instances, one can solve the isomorphism problem very efficiently using various heuristics. But a
polynomial-time algorithm working for all graphs is not known and it is very desirable to understand
the complexity of graph isomorphism. It is known that testing graph isomorphism is equivalent to
testing isomorphism of general mathematical structures~\cite{hedrlin}.  The notion of isomorphism is
widely used in mathematics when one wants to show that two seemingly different mathematical
structures are the same. One proceeds by guessing a mapping and proving that this mapping is an
isomorphism. The natural complexity question is whether there is a better way in which one
algorithmically derives an isomorphism. Similarly, regular covering is a well-known mathematical
concept which is algorithmically interesting and not understood.

Further, a regular covering is described by a semiregular subgroup of the automorphism group
$\Aut(G)$. Therefore it seems to be closely related to computation of $\Aut(G)$, since one should
have a good understanding of this group first, to solve the regular covering problem. The problem of
computing automorphism groups is known to be closely related to graph
isomorphism~\cite{mathon_isocount}. For survey of results about graph isomorphism,
see~\cite{babai1996automorphism,kkz}.

We note that Theorem~\ref{thm:algorithmic} allows to recognize finite planar Cayley graphs in
polynomial time. These graphs were already characterized by Maschke~\cite{maschke1896representation}
in 1896.  Unfortunately, finite planar Cayley graphs $G$ are very limited: either $G$ is $K_2$, or a
cycle, or a 3-connected planar graph. So $\Aut(G)$ is a \emph{spherical group} which is very simple.
Therefore, $G$ is either finite (with $\bv(G) \le 120$), representing one of the sporadic groups
(for instance, a truncated dodecahedron is a Cayley graph of $\gA_5$), or very simple (a cycle, a
prism, an antiprism, e.g.).  For study of infinite planar Cayley graphs,
see~\cite{droms1998connectivity} and the references therein.
 
\subsection{Structural Consequences} \label{sec:structural_cons}

\heading{A Proof of Negami Theorem.}
In 1988, Seiya Negami~\cite{negami} proved that a connected graph $H$ has a finite regular planar
cover $G$ if and only if $H$ is projective planar. If the graph $G$ is 3-connected, then $\Aut(G)$
is a \emph{spherical group}. Therefore the theorem can be easily proved using geometry, as we
discuss in Section~\ref{sec:planar_graphs}. The hard part of the proof is to deal with graphs $G$
containing 1-cuts and 2-cuts. Negami proves this by induction according to the size of $G$.  He
locates a minimal induced subgraph separated by 1-cut or 2-cut and replaces it in $G$, making a
smaller graph $G'$ from it. From induction hypothesis, every quotient $H'$ of $G'$ is planar or
projectively planar, and it is argued that a planar or projective planar embedding of the quotient
$H$ of $G$ can be produced from the embedding of $H'$.

While Negami's proof is quite short, the reader might ask several natural questions which are not
sufficiently answered:
\begin{packed_itemize}
\item While 3-connected graphs serve as starting points of the induction, what precisely is the role
of geometry and planarity when 1-cuts and 2-cuts appear?
\item Given a planar graph $G$, what is the full list of all possible quotients $H$ of $G$?
\item How is it determined whether a quotient $H$ is a planar or a projectively planar graph?
\end{packed_itemize}
Inspired by Negami~\cite{negami}, our work goes further and answers all these questions. Since
quotients of 3-connected planar graphs can be described geometrically (see
Section~\ref{sec:planar_graphs}), Theorem~\ref{thm:quotient_expansion} describes all quotients of
planar graphs.  Also, we work with a more general definition of regular graph coverings which admits
multigraphs with loops and half-edges. For instance, only half of the quotients of the cube in
Fig.~\ref{fig:quotients_of_cube} are admissible in Negami's definition. In the topological sense, we
allow non-isolated branch points which might also be placed in centers of edges. For a detailed
discussion, see Section~\ref{sec:quotients_of_planar_graphs}.

\heading{Characterizing Automorphism Groups of Planar Graphs.}
The key property that our reductions preserve essential information about symmetries of the graph
$G$ can be used to describe the automorphism groups of planar graphs.  This shows that computing
automorphism groups can be reduced to computing them for 3-connected graphs which are the spherical
groups.  The automorphism groups of planar graphs were non-inductively determined by
Babai~\cite{babai1975automorphism,babai1996automorphism} using a similar approach. An idea how to
compute a generating set of these groups was described by Colbourn and Booth~\cite{colbourn_booth},
but it was never fully developed. 

Our reduction is used in~\cite{knz} to obtain the first inductive characterization of the automorphism
groups of planar graphs. First, Proposition~\ref{prop:group_extension} is strengthened for planar
graphs to show that $\Aut(G_i) = \Aut(G_{i+1}) \ltimes \Ker(\Phi_i)$. After that the following
inductive characterization of stabilizers of vertices in connected planar graphs
$\Fix(\connectedplanar)$ is described. It is similar to Jordan-like characterization of the
automorphism groups of trees~\cite{jordan}. In what follows we shall use
some standard notation from group theory, see the next section for the exact definitions. 

\begin{theorem}[Klav\'{\i}k et al.~\cite{knz}] \label{thm:fix_groups}
The class $\Fix(\connectedplanar)$ is obtained inductively as follows:
\begin{packed_itemize}
\item[(a)] $\{1\} \in \Fix(\connectedplanar)$.
\item[(b)] If $\Psi_1,\Psi_2 \in \Fix(\connectedplanar)$, then $\Psi_1 \times \Psi_2 \in
\Fix(\connectedplanar)$.
\item[(c)] If $\Psi \in \Fix(\connectedplanar)$, then $\Psi \wr \gS_n, \Psi \wr \gC_n \in
\Fix(\connectedplanar)$.
\item[(d)] If $\Psi_1,\Psi_2,\Psi_3 \in \Fix(\connectedplanar)$, then
$$(\Psi_1^{2n} \times \Psi_2^n \times \Psi_3^n) \rtimes \gD_n \in \Fix(\connectedplanar),
\qquad\forall n \text{ odd}.$$
\item[(e)] If $\Psi_1,\Psi_2,\Psi_3,\Psi_4,\Psi_5 \in \Fix(\connectedplanar)$, then
$$(\Psi_1^{2n} \times \Psi_2^n \times \Psi_3^n \times \Psi_4^n \times \Psi_5^n) \rtimes \gD_n \in
\Fix(\connectedplanar),\qquad \forall n \ge 4, \text{ even}.$$ 
\item[(f)] If $\Psi_1,\Psi_2,\Psi_3,\Psi_4,\Psi_5,\Psi_6 \in \Fix(\connectedplanar)$, then
$$(\Psi_1^4 \times \Psi_2^2 \times \Psi_3^2 \times \Psi_4^2 \times \Psi_5^2 \times \Psi_6)
\rtimes \gC_2^2 \in \Fix(\connectedplanar).$$
\end{packed_itemize}
\end{theorem}

Next, the class of automorphism groups of connected planar graphs, denoted $\Aut(\connectedplanar)$,
are characterized as follows:

\begin{theorem}[Klav\'{\i}k et al.~\cite{knz}] \label{thm:planar_aut_groups}
Let $G$ be a planar graph with colored vertices and colored (possibly oriented) edges, which is
either 3-connected, or $K_1$, or $K_2$, or a cycle $C_n$. Let $m_1,\dots,m_\ell$ be the sizes of the
vertex- and edge-orbits of the action of $\Aut(G)$.  Then for all choices $\Psi_1,\dots,\Psi_\ell
\in \Fix(\connectedplanar)$, we have
$$(\Psi_1^{m_1} \times \cdots \times \Psi_\ell^{m_\ell}) \rtimes \Aut(G) \in
\Aut(\connectedplanar),$$
where the action of $\Aut(G)$ on the factors of each $\Psi_i^{m_i}$, $i = 1,\dots,\ell$, is induced by
the action of $\Aut(G)$ on the vertices and edges of $G$.

On the other hand, every group of $\Aut(\connectedplanar)$ can be constructed in the above way as
$$(\Psi_1^{m_1} \times \cdots \times \Psi_\ell^{m_\ell}) \rtimes \Sigma,$$
where $\Psi_1,\dots,\Psi_\ell \in \Fix(\connectedplanar)$ and $\Sigma$ is a spherical group.
\end{theorem}

This characterization leads to a quadratic-time algorithm for computing these automorphism groups.
The homomorphisms from the spherical groups $\cong\gD_n$ in operations (d), (e) and (f) in
Theorem~\ref{thm:fix_groups} and from $\Sigma$ in Theorem~\ref{thm:planar_aut_groups} determining
the semidirect products are induced by the action on the edges of the respective polyhedral graphs,
see \cite{knz} for details.

\section{Definitions and Preliminaries} \label{sec:preliminaries}

\subsection{Model of Graph} \label{sec:graph_model}

We first motivate the definition of extended graphs used in this paper.  The concept of graph
covering comes from topological graph theory where graphs are understood as 1-dimensional
CW-complexes. The main idea is that edges are represented by real open intervals and vertices by
points. The topological closure of an edge $e$ is either a closed interval, or a simple cycle.
In the first case, $e$ joins two different vertices $u$ and $v$ incident to $e$.  In the second
case, $e$ is incident just to one vertex $v$ and $e$ is a loop based at $v$.  When one considers
regular quotients of graphs, a third type of ``edges'' may appear~\cite{mns}.  For a non-trivial
involution swapping the end-vertices of an edge $e$, the regular covering projection $p$ maps $e$ to an
``edge'' $p(e)$ whose one end is incident to a vertex while the other is free. The topological
closure of $p(e)$ is homeomorphic to a half-closed interval which behaves as a ``half-edge''.

\begin{figure}[t!]
\centering
\includegraphics{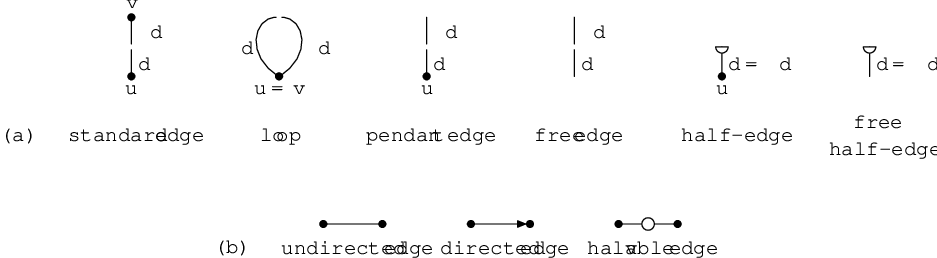}
\caption{(a) Four kinds of edges and two kinds of half-edges are depicted. We highlight two darts
composing each edge by a small gap, omitted in the remaining figures. To distinguish pendant edges
from half-edges, we end the latter by half-circles. (b) Three possible types for standard edges and
loops (pendant edges are always undirected). We note that only halvable edges may be projected to
half-edges which corresponds to cutting the middle circle in half, explaining the symbol for
half-edges.}
\label{fig:graph_definition}
\end{figure}

\heading{Definition of Extended Graphs.} An \emph{extended multigraph} $G$ (or just a \emph{graph}) is a
tuple $(\bD, \bV, \iota, \lambda)$, where $\bD$ is a set of darts, $\bV$ is a set of vertices,
$\iota : \bD \to \bV$ is a partial function of incidence, and $\lambda : \bD \to \bD$ is an
involution, pairing darts.  The set of edges $\bE$ is formed by orbits of $\lambda$ of size 2,
while orbits of size 1 form half-edges.  Each edge $\{d, \lambda d\}$ is one of the four
kinds:
\begin{packed_itemize}
\item a \emph{standard edge} if $\iota(d) \ne \iota(\lambda d)$,
\item a \emph{loop} if $\iota(d) = \iota(\lambda d)$,
\item a \emph{pendant edge} if exactly one of $\iota(d)$ and $\iota(\lambda d)$ is not defined, and
\item a \emph{free edge} if both $\iota(d)$ and $\iota(\lambda d)$ are not defined.
\end{packed_itemize}
For a half-edge $\{d\}$, we have $d = \lambda d$ and it is called a \emph{free half-edge} when
$\iota(d)=\iota(\lambda d)$ is not defined. See Fig.~\ref{fig:graph_definition}a.
\medskip

We introduce some further notation. A pendant edge attached to $v$ is called a \emph{single pendant
edge} if it is the only pendant edge attached to $v$. When we work with several graphs, we use
$\bD(G)$, $\bV(G)$, and $\bE(G)$ to denote the sets of darts, vertices and edges of $G$,
respectively. We denote $|\bD(G)|$ by $\bd(G)$, $|\bV(G)|$ by $\bv(G)$, $|\bE(G)|$ by $\be(G)$.
When $G$ contains no half-edges, clearly $\bd(G) = 2\be(G)$.  We consider graphs with colored
edges of three different edge types: \emph{directed}, \emph{undirected}, and a special type
called \emph{halvable}; see Fig.~\ref{fig:graph_definition}b.  It might seem strange to work with
such general objects. But when we apply reductions, we replace subgraphs by edges whose colors and
types encode isomorphism classes and symmetry types of replaced subgraphs. Even if $G$ and $H$ are
standard simple graphs, the more general colored multigraphs are naturally constructed in the
process of reductions.

Most graphs in this paper are assumed to be connected, so they contain no free edges and free
half-edges.  (Sometimes we consider subgraphs which may be disconnected and may contain them.) Note
that the standard concepts from graph theory such as vertex degree, connectedness, etc. easily
translate to extended multigraphs. 

When $A$ is a subset of $\bV(G) \cup \bE(G) \cup \bD(G)$, we denote by $G \setminus A$ the subgraph
created from $G$ by removing all elements of $A$ from $\bV(G) \cup \bE(G) \cup \bD(G)$, and by
modifying the incidence function $\iota$ and the pairing involution $\lambda$ in an obvious way.  In
particular, when we remove $u \in \bV(G)$, its incident edges are preserved in $G$ as pendant or
free edges. Similarly, when a dart $d \in \bD(G)$ is removed, the paired dart $\lambda d$, when $d
\ne \lambda d$, is preserved in $G$. When $A = \{x\}$, we often write $G \setminus x$ instead of $G
\setminus \{x\}$.

A vertex $u \in \bV(G)$, called an \emph{articulation}, forms a \emph{1-cut} $\{u\}$ in $G$ when
$G \setminus u$ is disconnected.  Similarly, two vertices $u,v \in \bV(G)$, $u \ne v$, form a
\emph{2-cut} $U = \{u,v\}$ in $G$ if $G \setminus U$ is disconnected. We say that a 2-cut $U$ is
\emph{non-trivial} if $\deg(u) \ge 3$ and $\deg(v) \ge 3$. An edge $e \in \bE(G)$ is called a
\emph{bridge-edge} when $G \setminus e$ is disconnected.

The induced subgraph consisting of two adjancent vertices joined by at least two edges is
called a \emph{dipole}.

For a subgraph $A$ of $G$, we use the topological notation to denote the \emph{boundary}
$\bo A$  and the \emph{interior} $\int A$ of $A$.  We set $\bo A$ equal to the set of vertices of
$A$ which are incident with an edge not contained in $A$. For the interior, we use the standard
topological definition $\int A = A \setminus \bo A$.

\subsection{Automorphism Groups} \label{sec:automorphisms}

\heading{Groups.}
For undefined concepts and results from permutation group theory, the reader is referred
to~\cite{rotman}.  We denote groups by Greek letters as for instance $\Psi$ or $\Gamma$. We use the
following notation for some standard families of groups:
\begin{packed_itemize}
\item $\gS_n$ for the symmetric group of all $n$-element permutations,
\item $\gC_n$ for the cyclic group of integers modulo $n$,
\item $\gD_n$ for the dihedral group of the symmetries of a regular $n$-gon, and
\item $\gA_n$ for the alternating group of all even $n$-element permutations.
\end{packed_itemize}

In this paper, by a group we usually mean a group of automorphisms of a graph acting on the set of
darts.  A group $\Psi$ \emph{acts} on a set $S$ in the following way. Each $g \in \Psi$
permutes the elements of $S$, and the \emph{action} is described by a mapping $\cdot : \Psi \times S
\to S$ where $1 \cdot x = x$ and $(gh) \cdot x = g \cdot (h \cdot x)$. 

\heading{Homomorphisms.}
We state the definitions in a very general setting of multigraphs with half-edges.  Let $G$ and $H$
be two graphs. Assuming that every vertex of $G$ has at least one dart incident, a homomorphism $h :
G \to H$ is fully described by a mapping $h_d : \bD(G) \to \bD(H)$ preserving edges and incidences
between darts and vertices, i.e.,
\begin{equation} \label{eq:homomorphism}
h_d(\lambda d) = \lambda h_d(d)\quad\text{and}\quad h_d(\iota(d)) = \iota(h_d(d)),\qquad \forall d
\in \bD(G)
\end{equation}
(where either both sides of the equations are defined, or none is defined). The mapping $h_d$
induces two mappings $h_v : \bV(G) \to \bV(G)$ and $h_e : \bE(G) \to \bE(G)$ connected together by
the very natural property $h_e(uv) = h_v(u) h_v(v)$ for every $uv \in \bE(G)$. (If some vertex $u$
has no dart incident, we also have to define $h_v(u)$.)

If $G$ is a simple graph (i.e., a graph containing only standard edges and no multiple edges are
allowed), then $h$ is determined by the action on the vertices, as is expected.  In most situations,
we omit subscripts and simply use $h(u)$, $h(d)$, or $h(uv)$. In addition, for colored graphs with
three edge types, we require that homomorphisms always preserves the colors, the edge types and the
direction of oriented edges.

For $A \subseteq \bV(G) \cup \bE(G) \cup \bD(G)$, we denote the \emph{restricted homomorphism} $A
\to H$ by $h |_A$.

\heading{Automorphism Groups.} For a graph $G$, an \emph{automorphism} $\pi$ is a homomorphism $G
\to G$ such that the mappings $\pi_d$, $\pi_v$, and $\pi_e$ are bijective.  The group of all
automorphisms of a graph $G$ will be denoted by $\Aut(G)$.  Each element $\pi \in \Aut(G)$  permutes
the vertices, edges and darts such that the edges and the incidences between the darts and the
vertices are preserved.

For $x \in \bV(G) \cup \bE(G) \cup \bD(G)$, the \emph{orbit} $[x]$ in the action of $\Psi \le
\Aut(G)$ is the set $\{\pi(x) \mid \pi \in \Psi\}$, and $[x]$ is called a vertex-, an edge-, or a
dart-orbit if $x$ is a vertex, an edge, or a dart, respectively.  The \emph{stabilizer} of $x \in
\bV(G) \cup \bE(G) \cup \bD(G)$ in the action of $\Psi \le \Aut(G)$ is the subgroup $\{\pi \mid
\pi(x) = x\}$, and again, it is called a vertex-, an edge-, or a dart-stabilizer if $x$ is a vertex,
an edge, or a dart, respectively.  An action of $\Psi\leq \Aut(G)$ is called \emph{semiregular} if
it has no non-trivial (i.e., non-identity) dart- and vertex-stabilizers.  By definition, an
edge-stabilizer in semiregular action is either trivial or isomorphic to $\gC_2$. We require an
edge-stabilizer in a semiregular action to be trivial, unless it is a halvable edge, when the action
may contain an involution transposing the two darts. We say that a group is \emph{semiregular} if
its action is semiregular. Through the paper, the letter $\Gamma$ is reserved for semiregular
subgroups of $\Aut(G)$. We say that $\pi \in \Aut(G)$ is semiregular if the subgroup
$\left<\pi\right>$ is semiregular. (Note that this is equivalent to the fact that $\pi$ has all its
cycles of the same length.)

For a set $S \subseteq \bV(G) \cup \bE(G) \cup \bD(G)$, the \emph{point-wise stabilizer} of $S$ in
$\Psi \le \Aut(G)$ is the subgroup of $\Psi$ consisting of all automorphisms $\pi$ such that $\pi(x)
= x$ for all $x \in S$, while the \emph{set-wise stabilizer} of $S$ in $\Psi \le \Aut(G)$ consists
of all automorphisms $\pi$ such that $\pi(x) \in S$ for all $x \in S$.  When working with subgraphs
$A$ of $G$, we consider the following two boundary-preserving subgroups of $\Aut(A)$:
\begin{packed_itemize}
\item $\Fix(\bo A)$ is the point-wise stabilizer of $\bo A$ in $\Aut(A)$.
\item $\Autbo A$ is the set-wise stabilizer of $\bo A$ in $\Aut(A)$.
\end{packed_itemize}
Observe that $\Fix(\bo A)$ considered as a subgroup of $\Aut(G)$ is equivalently the point-wise
stabilizer of $G \setminus \int A$ in $\Aut(G)$.

\heading{Isomorphisms.}
For graphs $G$ and $H$, an \emph{isomorphism} $\sigma : G \to H$ is a bijective homomorphism.
Observe that an isomorphism $G \to G$ is an automorphism of $G$. We denote existence of an
isomorphism between $G$ and $H$ by $G \cong H$. Naturally, for colored graphs with three edge types,
isomorphisms are required to preserve the colors, the edge types and the directions of oriented
edges. 

For subgraphs of $G$, we usually consider only isomorphisms preserving their boundaries. Let $A$,
$A'$ be subgraphs of $G$. An isomorphism $\sigma : A \to A'$ is called a \emph{$\bo$-isomorphism} if
$\sigma(\bo A) = \bo \sigma(A)$. If such a $\bo$-isomorphism exists, we say that $A$ is
$\bo$-isomorphic to $A'$, denoted $A \cong_\bo A'$. Observe that for every subgraph $A$ and every
automorphism $\pi \in \Aut(G)$, the restriction $\pi |_A$ is a $\bo$-isomorphism from $A$ to $\pi(A)$.

\begin{figure}[t!]
\centering
\includegraphics{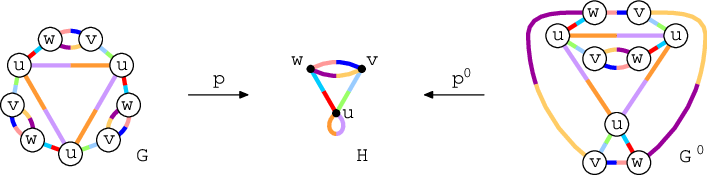}
\caption{Graphs $G$ and $G'$ covering a graph $H$ with covering projections $p$ and $p'$,
respectively. The projections $p_v$ and $p'_v$ take a vertex with a label $x$ onto a vertex with the
same label $x\in \bV(H)$ and the projections $p_d$ and $p'_d$ take each dart $\bD(G)$ and $\bD(G')$
onto a dart of $\bD(H)$ with the same color.  Notice that the loop in $H$ lifts along $p$ and $p'$
into three standard edges joining the vertices with the same label, while the two parallel edges
joining $v$ to $w$ in $H$ lift along $p$ to a $2$-factor consisting of three parallel edges in $G$,
and they lift along $p'$ to $2$-factor consisting of a parallel edge and a $4$-cycle.}
\label{fig:cover_examples}
\end{figure}

\begin{wrapfigure}{r}{0.5\textwidth}
\begin{center}
\vskip -4ex
\includegraphics{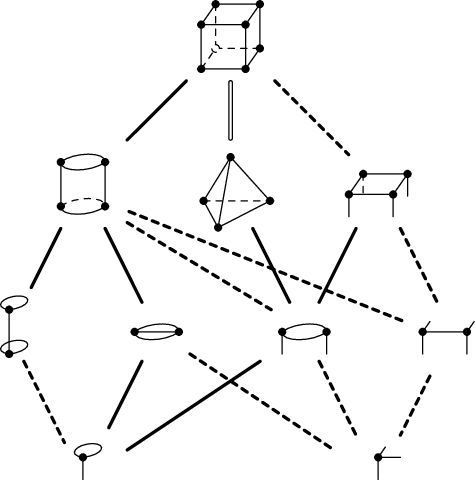}
\vskip -2ex
\end{center}
\caption{The Hasse diagram of all quotients of the cube graph depicted in a geometric way. When
semiregular actions fix edges, the quotients contain half-edges. The quotients connected by bold
edges are obtained by 180 degree rotations. The quotients connected by dashed edges are obtained by
reflections. The tetrahedron is obtained by the antipodal symmetry of the cube, and its
half-quotient is obtained by a 180 degree rotation with the axis going through the centers of two
non-incident edges of the tetrahedron.}
\label{fig:quotients_of_cube}
\vskip -5ex
\end{wrapfigure}

\subsection{Coverings}

A graph $G$ \emph{covers} a graph $H$ (or $G$ is a \emph{cover} of $H$) if there exists a locally
bijective homomorphism $p$ called a \emph{covering projection}. The property to be locally bijective
states that for every vertex $u \in \bV(G)$ the mapping $p_d$ restricted to the darts incident with
$u$ is a bijection.  Figure~\ref{fig:cover_examples} contains two examples of graph covers. Observe
that each vertex, dart, or an edge has exactly three preimages.

\heading{Fibers.}
A \emph{fiber} over a vertex $v \in \bV(H)$ is the set $p^{-1}(v)$, i.e., the set of all vertices
$\bV(G)$ that are mapped to $v$, and similarly for fibers over edges and darts.  We adopt the
standard assumption that both $G$ and $H$ are connected.  It is well known that all fibers of $p$ are of
the same size.  In other words, $\bd(G) = k \cdot \bd(H)$ and $\bv(G) = k \cdot \bv(H)$ for some $k
\in \mathbb N$ which is the size of each fiber, and we say that $G$ is a \emph{$k$-fold cover} of
$H$.

\heading{Regular Coverings.} We aim to consider regular graph coverings which are closely related to
semiregular groups of automorphisms. 

Let $\Gamma$ be any semiregular subgroup of $\Aut(G)$. It defines a graph $G / \Gamma$ called a
\emph{regular quotient} (or simply \emph{quotient}) of $G$ as follows: The vertices of $G / \Gamma$
are the orbits of the action of $\Gamma$ on $\bV(G)$, the darts of $G / \Gamma$ are the orbits of
the action of $\Gamma$ on $\bD(G)$.  A vertex-orbit $[v]$ is incident with a dart-orbit $[d]$ if and
only if the vertices of $[v]$ are incident with the darts of $[d]$. (Because the action of
$\Gamma$ is semiregular, each vertex of $[v]$ is incident with exactly one dart of $[d]$, so
this is well defined.) We say that $G$ \emph{regularly covers} $H$ if there exists a regular
quotient of $G$ isomorphic to $H$. 

We naturally construct the \emph{regular covering projection} $p : G \to G / \Gamma$ by mapping the
vertices to its vertex-orbits and darts to its dart-orbits. Concerning an edge $e \in
\bE(G)$, it is mapped to an edge $p(e)$ of $G / \Gamma$ if the two darts of $e$ belong to different
dart-orbits of $\Gamma$. If they belong to the same dart-orbit, then $p(e)$ is a half-edge of $G / \Gamma$.
The projection $p$ is a $|\Gamma|$-fold regular covering projection.

From the two examples from Fig.~\ref{fig:cover_examples}, the covering projection $p$ is a $3$-fold
regular covering projection over $H$ while the covering projection $p'$ is an irregular $3$-fold
covering projection. Observe that $p$ is induced by a semiregular group $\Gamma \le \Aut(G)$ of
order three rotating the central triangle of $G$, while $G'$ does not admit a semiregular group
$\Gamma \le \Aut(G')$ of order three. Figure~\ref{fig:quotients_of_cube} depicts all the regular
quotients of the cube graph.

\begin{figure}[b!]
\centering
\includegraphics{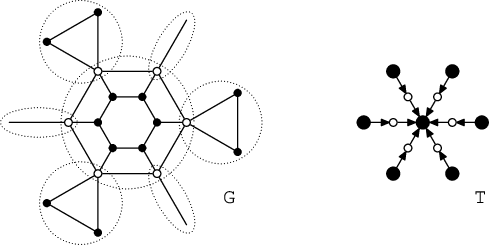}
\caption{On the left, an example graph $G$ with denoted blocks. On the right, the corresponding
block-tree $T$ is depicted, rooted at the central block. The white vertices correspond to the
articulations and the big black vertices correspond to the blocks.}
\label{fig:example_of_block_tree}
\end{figure}

\subsection{Block-trees and Their Automorphisms} \label{sec:block_trees}

The \emph{block-tree} $T$ of $G$ is a tree defined as follows. Consider all articulations in $G$ and
all maximal $2$-connected subgraphs which we call \emph{blocks} (with bridge-edges and pendant edges
also counted as blocks).  The block-tree $T$ is the incidence graph between articulations and
blocks. For an example, see Fig.~\ref{fig:example_of_block_tree}.  It is well known that every
automorphism $\pi \in \Aut(G)$ induces an automorphism $\pi' \in \Aut(T)$.

\heading{The Central Block.} For a tree, its \emph{center} is either the central vertex, or the
central pair of vertices of a longest path, depending on the parity of its length.  For the
block-tree $T$, all leaves are blocks and each longest path is of an even length.  Therefore, $T$
has a central vertex which is either a \emph{central articulation}, or a \emph{central block} of
$G$.

\begin{lemma} \label{lem:central_block}
If $G$ has a non-trivial semiregular automorphism, then $G$ has a central block.
\end{lemma}

\begin{proof}
For contradiction, suppose that $G$ has a central articulation $u$. Let $T$ be the block-tree. Every
automorphism of a tree preserves its center, so $\Aut(T)$ preserves $u$.  Also, all automorphisms of
$\Aut(G)$ preserve $u$ since every automorphism of $\Aut(G)$ induces an automorphism of $\Aut(T)$.
In particular, we get that a non-trivial semiregular automorphism $\pi$ fixes the vertex $u$. By
the semiregularity, $\pi$ is the identity which is a contradiction.\qed
\end{proof}

We orient the edges of the block tree $T$ towards the central vertex, so the block tree becomes
rooted.  A \emph{rooted subtree} $T'$ of the block tree $T$ is the induced subgraph of $T$
determined by a vertex $v$, called the \emph{root} of $T'$, and by all its predecessors.

Suppose that $G$ has a central block $C$.  Let $u$ be an articulation contained in $C$. By $T_u$, we
denote the rooted subtree of $T$ determined by $u$, and let $G_u$ be the graph induced by all
vertices of the blocks of $T_u$.

\begin{lemma} \label{lem:semiregular_action_on_blocks}
Let $\Gamma$ be a semiregular subgroup of $\Aut(G)$. If $u$ and $v$ are two articulations of the
central block and belonging to the same orbit of $\Gamma$, then $G_u \cong_\bo G_v$. Moreover there is a
unique $\pi \in \Gamma$ which maps $G_u$ to $G_v$ and $\bo G_u$ to $\bo G_v$.
\end{lemma}

\begin{proof}
Notice that either $G_u = G_v$, or $G_u \cap G_v = \emptyset$. Since $u$ and $v$ are in the same
orbit of $\Gamma$, there exists $\pi \in \Gamma$ such that $\pi(u) = v$. Consequently $\pi(G_u)
= G_v$.  Suppose that there exist $\pi,\sigma \in \Gamma$ such that $\pi(G_u) = \sigma(G_u) = G_v$.
Then $\pi\cdot\sigma^{-1}$ is an automorphism of $\Gamma$ fixing $u$. Since $\Gamma$ is semiregular,
$\pi = \sigma$.\qed
\end{proof}

In the language of quotients, it means that $G / \Gamma$ consists of the quotient $C / \Gamma$ of
the central block, together with the graphs $G_u$ attached to $C / \Gamma$, one copy for each orbit
$[u]$ of  $\Gamma$, where $u$ ranges through all articulations in $C$.

\heading{Why Not Just 2-connected Graphs?}
Since the behaviour of regular coverings with respect to 1-cuts in $G$ is very simple, a natural
question follows: why do we not restrict ourselves to 2-connected graphs $G$?  The issue is that the
quotient $C / \Gamma$ might not be 2-connected (see Fig.~\ref{fig:atom_projection_cases} on the
right), so it may consists of many blocks in $H$. When $H$ contains a rooted subtree of blocks
isomorphic to $G_u$, it may correspond to $G_u$, or it may correspond to a quotient of a subgraph $C
/ \Gamma$, together with some other $G_v$ attached.  Therefore, we work with 1-cuts together with
2-cuts and we define the 3-connected reduction for 1-cuts in $G$ as well, unlike 
in~\cite{maclane,trakhtenbrot,tutte_connectivity,quadratic_isomorphism_planar,hopcroft_tarjan_dividing,cunnigham_edmonds,walsh,bienstock,
droms1995structure,spqr1,spqr2,spqr3,spqr_linear}.  This modification is essential for the algorithm
for regular covering testing described in~\cite{fkkn,fkkn16}.

Further, we believe this unified approach for 1-cuts and 2-cuts is preferable in general.
Our definition of atoms and reductions is slightly more complicated. But we get one
reduction procedure and one reduction tree for the entire connected graph $G$. This makes the
characterization of automorphism groups of planar graphs in~\cite{knz} more understandable.

\section{Structural Properties of Atoms} \label{sec:atoms}

In this section, we introduce special inclusion-minimal subgraphs of $G$ called atoms. We
investigate their structural properties, in particular their behaviour with respect to regular
covering projections.

\subsection{Definition of Atoms} \label{sec:definition_of_atoms}

\begin{figure}[b!]
\centering
\includegraphics{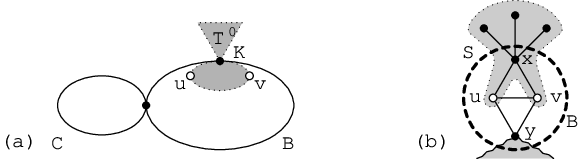}
\caption[]{(a) In the definition of a proper part, notice that a 2-cut $U = \{u,v\}$ is defined with
respect to some block $B$, but $K$ (depicted in gray) is defined with respect to $G$ as a connected
component of $G \setminus U$. So $K$ is a connected component of $B \setminus U$ together with the
subgraphs induced by some attached rooted subtrees. Assuming that the graph has a central block $C$,
we require that $C$ does not belong to $K$, so either $B = C$, or $C$ is some other block outside of
$K$.\\
(b) A section of a graph $G$, with the central block/articulation separated by $y$. While $X =
\{x,y\}$ forms a non-trivial 2-cut in $G$ (both $x$ and $y$ are articulations), $X$ is not used in
the definition of proper parts since $X$ is not a 2-cut within the block $B$. The only non-trivial
2-cut within the block $B$ is $U = \{u,v\}$. There are three connected components in $G \setminus U$
and only the one containing $x$ forms together with $U$ the proper part $S$, highlighted in gray.
The free edge is forbidden by the definition and the connected component of $G \setminus U$ contains
the central block/articulation, which is also forbidden.}
\label{fig:proper_parts}
\end{figure}

\begin{figure}[t!]
\centering
\includegraphics{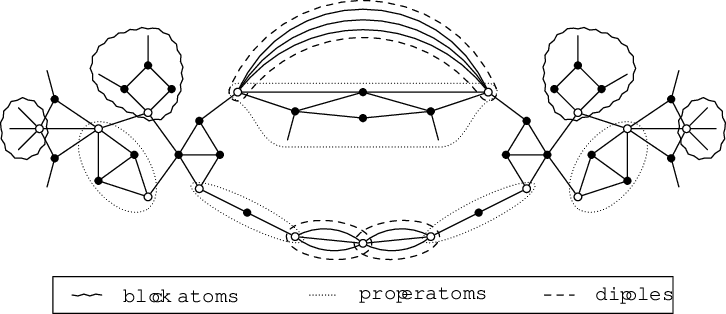}
\caption{An example of a graph with denoted atoms. The white vertices belong to the boundary of some
atom, possibly several of them.}
\label{fig:atoms_examples}
\end{figure}

We first define a set $\calP$ of subgraphs of $G$ called \emph{parts} which are candidates for
atoms:
\begin{packed_itemize}
\item A \emph{block part} is a subgraph $S$ of $G$ defined by a rooted subtree $T'$ determined by a
vertex $x$ of the block-tree other than the central block. The subgraph $S$ consists of all blocks
in $T'$ and we require that $S$ is non-isomorphic to a single pendant edge.
\item A \emph{proper part} is a subgraph $S$ of $G$ defined by a non-trivial 2-cut $U = \{u,v\}$ in
some block $B$ of $G$. The subgraph $S$ consists of a connected component $K$ of $G \setminus U$
together with $u$ and $v$. In addition, we require that $K$ is not a single free edge and $K$ does
not contain all vertices of the central block, or the central articulation.
Figure~\ref{fig:proper_parts} illustrates this definition.
\item A \emph{dipole part} is any dipole.
\end{packed_itemize}

\begin{figure}[b!]
\centering
\includegraphics{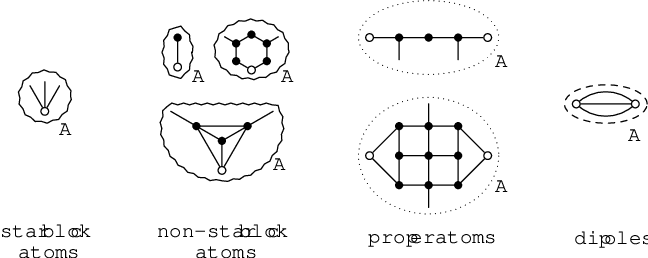}
\caption{Examples of four different types of atoms, with white vertices belonging to $\bo A$.}
\label{fig:types_of_atoms}
\end{figure}

\noindent The inclusion-minimal elements of $\calP$ are called \emph{atoms}. We distinguish
\emph{block atoms}, \emph{proper atoms} and \emph{dipoles} according to the type of the defining
part.  Block atoms are either pendant stars called \emph{star block atoms}, or pendant blocks
possibly with single pendant edges attached to them called \emph{non-star block atoms}.
Figure~\ref{fig:atoms_examples} gives an example and Figure~\ref{fig:types_of_atoms} an
overview of different types of atoms. 

The above concepts of a proper atom and dipoles have their counter-parts in the literature, they are
called pseudo-bricks and bonds, respectively~\cite{walsh}.  Some of the following properties and
results can be found in literature,
see~\cite{trakhtenbrot,tutte_connectivity,hopcroft_tarjan_dividing,cunnigham_edmonds,walsh} for
instance.  The novelty of our approach is the use of pendant edges which allow to define atoms also
for 1-connected graphs. For the reader's convenience, we prove all the properties of the atoms used
in further considerations to make this paper self-contained.

Notice that each proper atom is a subgraph of some block, together with some single pendant edges
attached to it. A dipole part is by definition always inclusion-minimal, and therefore it is an
atom. We have $|\bo A| = 1$ for a block atom $A$, and $|\bo A| = 2$ for a proper atom or dipole $A$.
The interior of a star block atom or a dipole is a set of free edges. Observe that for a proper atom
$A$,  the vertices of $\bo A$ are exactly the vertices $\{u,v\}$ of the non-trivial 2-cut $U$ used
in the definition of proper parts. Also the vertices of $\bo A$ of a proper atom are never adjacent
in $A$ (while they might be adjacent in $G$) since otherwise $K = \int A$ would be disconnected or
it would consist of only one free edge. Further, no block or proper atom contains parallel edges;
otherwise a dipole would be its subgraph, so it would not be inclusion minimal.

It is important to point out why single pendant edges play a special role in the definition of
parts. In Section~\ref{sec:reduction_and_expansion}, we define the reduction process which replaces
block atoms by pendant edges. We use single pendant edges since multiple non-star block atoms may
share the articulation in their boundaries, resulting into a star of pendant edges in the reduced
graph which more fatefully preserves the structure of the original graph. For a detailed
discussion, see~\cite[Section 7.3.2]{phd_thesis}.

\subsection{Primitive Graphs} \label{sec:primitive_graphs}

A graph is called \emph{primitive} if it contains no atoms.  The following lemma characterizing
primitive graphs can be alternatively obtained from the well-known theorem by
Trakhtenbrot~\cite{trakhtenbrot};\footnote{We consider $K_1$ with an attached single pendant edge
as a graph with a central articulation.}
see Fig.~\ref{fig:primitive_graphs}.

\begin{lemma} \label{lem:primitive_graphs}
Let $G$ be a primitive graph. If $G$ has a central block, then it is a 3-connected
graph, a cycle $C_n$ for $n \ge 2$, or $K_2$, or can be obtained from the aforementioned graphs by attaching
single pendant edges to at least two vertices. If $G$ has a central articulation, then it is $K_1$,
possible with a single pendant edge attached.
\end{lemma}

\begin{proof}
Clearly, the graphs mentioned in the statement are primitive.  On the other hand, a primitive graph
$G$ has a central block/articulation. All blocks attached to it have to be single pendant edges,
otherwise $G$ would contain a block atom. If $G$ has a central articulation $u$, then all the
connectivity components separated by $u$ are free edges.  In particular, all the edges are pendant,
and since $G$ is primitive, the number of pendant edges is at most one.  Thus $G$ is $K_1$, possibly
with a single pendant edge attached.  If $G$ has a central block, after removing all pendant edges,
we get the 2-connected graph $B$ consisting of only the central block. We argue that $B$ is one of
the stated graphs.

\begin{figure}[t!]
\centering
\includegraphics{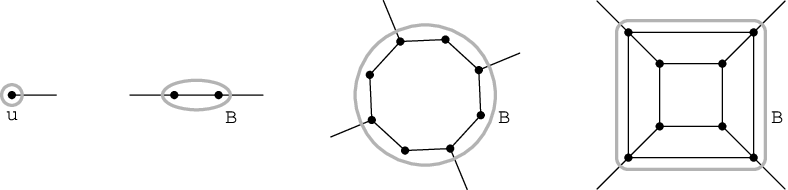}
\caption{A primitive graph with a central articulation is $K_1$, and with a central block is either
$K_2$, $C_n$, or a 3-connected graph, in all these cases with possible single pendant edges attached
to it.}
\label{fig:primitive_graphs}
\end{figure}

Now, let $u$ be a vertex of the minimum degree in $B$. If $\deg(u)=1$, the graph $B$ has to be
$K_2$, otherwise it would not be 2-connected. If $\deg(u)=2$, then either the graph $B$ is a cycle
$C_n$, or $u$ is an inner vertex of a path connecting two vertices $x$ and $y$ of degree at least
three such that all inner vertices are of degree two. But then this path is a proper atom which is a
contradiction. Finally, if $\deg(u) \ge 3$, then every 2-cut is non-trivial, and since $B$ contains
no proper atoms, $B$ has to be 3-connected.\qed
\end{proof}

\subsection{Structure of Atoms} \label{sec:structure_of_atoms}

We call a graph \emph{essentially 3-connected} if it is a 3-connected graph possibly with some
single pendant edges attached to it.  Similarly, a graph is called \emph{essentially a cycle} if it
is a cycle possibly with some single pendant edges attached to it.  The structure of star block
atoms and dipoles is already clear from the definition.  The following lemmas describe the structure
of non-star block and proper atoms; different examples are depicted in Fig.~\ref{fig:types_of_atoms}.

\begin{lemma} \label{lem:non_star_block_atoms}
Every non-star block atom $A$ is either $K_2$, possibly with an attached single pendant edge, essentially a
cycle, or essentially 3-connected.
\end{lemma}

\begin{proof}
Clearly, the described graphs are possible non-star block atoms.  Since $A$ does not contain any
smaller block atom, then $A$ is $2$-connected graph, possibly with some single pendant edges
attached.  By removing all single pendant edges, we get a 2-connected graph $B$, otherwise $A$
contains a smaller block part, which is a smaller block part in $G$ as well.

Let $u$ be a vertex of the minimum degree in $B$.  We have $\deg(u) > 0$, otherwise $B = K_1$ and
$A$ is a pendant edge. If $\deg(u)=1$, the graph $B$ has to be $K_2$, otherwise it would not be
2-connected. If $\deg(u)=2$, then either the graph $B$ is a cycle $C_n$, or $u$ is an inner vertex
of a path connecting two vertices $x$ and $y$ of degree at least three such that all inner vertices
are of degree two. But then this path determines a proper atom in $B$ which is also a proper atom in
$G$, a contradiction. Finally, if $\deg(u) \ge 3$, then every 2-cut is non-trivial, and since $B$
contains no proper atoms, it has to be 3-connected.\qed
\end{proof}

Let $A$ be a proper atom with $\bo A = \{u,v\}$. We define the \emph{extended proper atom} $A^+$ as
$A$ with the additional edge $uv$. 

\begin{lemma} \label{lem:proper_atoms}
For every proper atom $A$, the extended proper atom $A^+$ is either essentially a cycle, or
essentially 3-connected.
\end{lemma}

\begin{proof}
Clearly, the described graphs are possible extended proper atoms $A^+$.  Notice that $A^+$ consists
of a $2$-connected graph, possibly with single pendant edges attached, otherwise $A$ contains a smaller
block part. By removing all single pendant edges, we get a 2-connected graph $B^+$, otherwise $A^+$
contains a smaller block part. Let $\bo A = \{u,v\}$, we have $\deg(u) \ge 2$ and $\deg(v) \ge 2$ in
$A^+$ (and their degrees are preserved in $B^+$).

Let $w$ be a vertex of the minimum degree in $B^+$.  We have $\deg(w) > 1$, otherwise $A$ again
contains a smaller block part.  If $\deg(w)=2$, then either the graph $B^+$ is a cycle $C_n$, or $w$
is an inner vertex of a path connecting two vertices $x$ and $y$ of degree at least three such that
all inner vertices are of degree two. But then this path is a proper atom in $A^+$. It corresponds
to a proper atom in the original graph since the edge $uv$ in $A^+$ corresponds to some path in $G
\setminus \int A$, so we get a contradiction with the minimality of $A$. Finally, if $\deg(w) \ge
3$, then every 2-cut is non-trivial, and since $B^+$ contains no atoms, it has to be
3-connected.\qed
\end{proof}

Notice that for all atoms $A$, only vertices of $\int A$ might have single pendant edges attached
for the following reason. Dipoles and star block atoms contain no single pendant edges from the
definition.  A proper atom $A$ is defined as a connected component $K$ of $G \setminus U$
together with $U$, for some non-trivial 2-cut $U$, and single pendant edges attached to $U$ do
not belong to $K$. A non-star block atom $A$ is the subgraph induced by a rooted subtree
determined by some block $B$, and a single pendant edges attached to $\bo A$ would be a sibling
of $B$ in the block tree, so it does not belong to $A$.

\begin{lemma} \label{lem:aut_ess_3-conn}
Let $A$ be an essentially 3-connected graph, and let $B$ arise from $A$ by removal of all the single
pendant edges of $A$. Then every automorphism $\Aut(A)$ is isomorphic to a subgroup of $\Aut(B)$.
Further, for every $S \subseteq \bV(B) \cup \bE(B) \cup \bD(B)$, the stabilizer of $S$ in $\Aut(A)$
is a subgroup of the stabilizer of $S$ in $\Aut(B)$.
\end{lemma}

\begin{proof}
These single pendant edges behave like markers at the vertices they are attached to, thus forming a
2-partition of $\bV(A)$ preserved by the every automorphism in $\Aut(A)$. For every
automorphism $\pi \in \Aut(A)$, the restriction $\pi |_B$ belongs to $\Aut(B)$, and either both, or
neither stabilizes $S$.\qed
\end{proof}

In particular, when $A$ is an atom, we get that $\Autbo{A} \le \Autbo{B}$ and $\Fix(\bo A) \le
\Fix(\bo B)$.

\subsection{Non-overlapping Atoms} \label{sec:non_overlapping_atoms}

Our goal is to replace atoms by edges, and so it is important to know that the atoms cannot overlap
too much. The reader can see  in Fig.~\ref{fig:atoms_examples} that the atoms only share their
boundaries. This is true in general, and we are going to prove it in two steps. 

\begin{figure}[b!]
\centering
\includegraphics{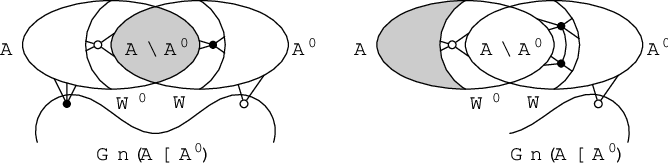}
\caption{We depict the vertices of $\bo A$ in black and the vertices of $\bo A'$ in white.
In both cases, we find a subset of $A$ belonging to $\calP$ (its interior is highlighted in gray).}
\label{fig:intersecting_interiors}
\end{figure}

\begin{lemma} \label{lem:nonintersecting_interiors}
The interiors of distinct atoms are disjoint.
\end{lemma}

\begin{proof}
For contradiction, let $A$ and $A'$ be two distinct atoms with non-empty intersections of $\int A$
and $\int A'$. First suppose that one of them, say $A$, is a block atom. Then $A$ corresponds to a
subtree of the block-tree which is attached by an articulation $u$ to the rest of the graph. If $A'$
is a block atom then it corresponds to some subtree, and we can derive that $A \subseteq A'$ or $A'
\subseteq A$. If $A'$ is a dipole, then it is a subgraph of a block, and thus subgraph of $A$.
If $A'$ is a proper atom, it is defined with respect to some block $B$. If $B$ belongs to the
subtree corresponding to $A$, then $A' \subseteq A$. Otherwise, a subtree of blocks containing $A$
is attached to $A'$, so $A \subseteq A'$. In both cases, we get a contradiction with the minimality.
Similarly, if one of the atoms is a dipole, we can easily derive a contradiction with the minimality.

The last case to consider is that both $A$ and $A'$ are proper atoms. Since the interiors are
connected and the boundaries are defined as neighbors of the interiors, it follows that both $W' = A
\cap \bo A'$ and $W = A' \cap \bo A$ are nonempty. We have two cases according to the sizes of these
intersections depicted in Fig.~\ref{fig:intersecting_interiors}.

If $|W| = |W'| = 1$, then $W \cup W'$ is a 2-cut separating $\int A \cap \int A'$ which contradicts
the minimality of $A$ and $A'$. Assume, without loss of generality, 
that $|W| = 2$. Then there is no edge
between $\int A \setminus (\int A' \cup W')$ and the remainder of the graph $G \setminus (\int A
\cup \int A')$. Therefore, $\int A \setminus (\int A' \cup W')$ is separated by a 2-cut $W'$ which
again contradicts the minimality of $A$. We note that in both cases the constructed 2-cut is
non-trivial since it is formed by vertices of non-trivial 2-cuts $\bo A$ and $\bo A'$.\qed
\end{proof}

Next we show a stronger version of the previous lemma which states that two atoms can intersect only
in their boundaries.

\begin{lemma} \label{lem:nonintersecting_atoms}
Let $A$ and $A'$ be two atoms. Then $A \cap A' = \bo A \cap \bo A'$.
\end{lemma}

\begin{figure}[t!]
\centering
\includegraphics{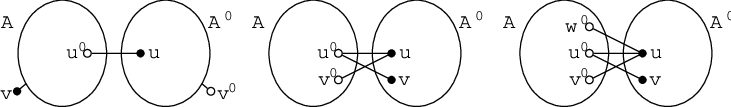}
\caption{An illustration of the main steps of the proof of Lemma~\ref{lem:nonintersecting_atoms}.
We depict the vertices of $\bo A$ in black and the vertices of $\bo A'$ in white.}
\label{fig:intersecting_atoms}
\end{figure}

\begin{proof}
We already know from Lemma~\ref{lem:nonintersecting_interiors} that $\int A \cap \int A'=
\emptyset$. It remains to argue that, say, $\int A \cap \bo A' = \emptyset$. If $A'$ is a block atom,
then $\bo A'$ is the articulation separating $A$. If $A$ contains this articulation in its interior,
it also contains $A'$ in its interior as well, contradicting $\int A \cap \int A' = \emptyset$. Similarly, if
$A$ is a block atom, then $A'$ has to be contained in $\int A$ or vice versa which again
contradicts $\int A \cap \int A' = \emptyset$.

It remains to argue the case when both $A$ and $A'$ are proper atoms or dipoles. Let $\bo A =
\{u,v\}$ and $\bo A' = \{u',v'\}$.  First we deal with the dipoles. If $A$ is a dipole, then $\int A$
contains no vertices, and consequently, the statement holds true.  If $A'$ is a dipole and $A$ is a
proper atom with $u' \in \int A$, then  the edges of $A'$ belong to $A$. Thus $A' \subsetneq A$,
contradicting the minimality.

We conclude with the remaining case assuming that both $A$ and $A'$ are proper atoms. Recall that
$\bo A$ is defined as the set neighbors of vertices of $\int A$ in $G\setminus \int A$, and
similarly, $\bo A'$ is the set of neighbors of vertices of $\int A'$ in $G\setminus \int A'$, for an
illustration see Fig.~\ref{fig:intersecting_atoms}.

Suppose for contradiction that $\int A \cap \bo A' \ne \emptyset$ and let $u' \in \int A$. By
definition, $u'$ has at least one neighbor in $\int A'$, and since $\int A \cap \int A' =
\emptyset$, this neighbor does not belong to $\int A$.
Therefore, without loss of generality, we have $u \in \int A'$ and $uu' \in \bE(G)$.
Since $A$ is a proper atom,  the set $\{u',v\}$ is not a 2-cut, so there is another
neighbor of $u$ in $\int A$, which has to be equal $v'$. Symmetrically, $u'$ has another neighbor in
$\int A'$ which is $v$.  So $\bo A \subseteq \int A'$ and $\bo A' \subseteq \int A$.  If $\bo A =
\int A'$ and $\bo A' = \int A$, the graph is $K_4$ (since the minimal degree of cut-vertices is
three) which contradicts existence of 2-cuts and atoms. If for example $\int A \ne \bo A'$, then
$\bo A'$ does not cut a subset of $\int A$, so there exists $w' \in \int A$ which is a neighbor of
$\int A'$, which contradicts that $\bo A'$ cuts $\int A'$ from the rest of the graph.\qed
\end{proof}

\subsection{Symmetry Types of Atoms}

\begin{figure}[b!]
\centering
\includegraphics{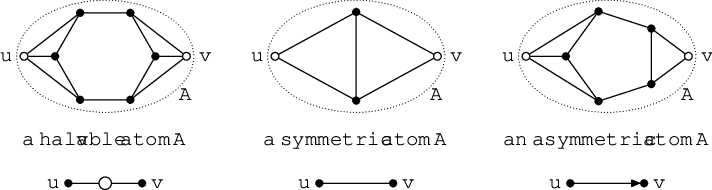}
\caption{The three types of atoms and the corresponding edge types which we use in the reduction.
White vertices belong to $\bo A$.}
\label{fig:symmetry_types_of_atoms}
\end{figure}

We distinguish three symmetry types of atoms. Recall that
$\Autbo A$ denotes the setwise stabilizer of $\bo A$ in $\Aut(A)$. If $A$ is a block atom, then its
symmetry type is by definition \emph{symmetric}. Let $A$ be a proper atom or a dipole with $\bo A = \{u,v\}$. Then we
distinguish the following three symmetry types, see Fig.~\ref{fig:symmetry_types_of_atoms}:
\begin{packed_itemize}
\item \emph{The atom $A$ is halvable} when there exits a semiregular involutory automorphism $\tau \in \Autbo A$
which exchanges $u$ and $v$.
\item \emph{The atom $A$ is symmetric} when it is not halvable, but there exists an automorphism in
$\Autbo A$ which exchanges $u$ and $v$.
\item \emph{The atom $A$ is asymmetric} when it is neither halvable, nor symmetric.
\end{packed_itemize}
We note that a symmetric proper atom $A$ might not have any involutory automorphism in $\Autbo A$
exchanging the boundary vertices; see~\cite{knz} for an example and further discussions.

When an atom is reduced, we replace it by an edge carrying the type. Therefore we work with
multigraphs with three edge types: \emph{halvable edges}, \emph{undirected edges} and \emph{directed
edges}.  Since our primary aim is to investigate action of the automorphism group of a graph, we
choose the orientation of the directed edges, replacing asymmetric atoms forming an orbit of the
automorphism group, to be compatible with the action of the group. 

For these multigraphs, we naturally consider only the automorphisms which preserve these edge types
and of course the orientation of directed edges. This generalized definition is used to define
symmetry types of their atoms. In the definition of a halvable atom, the automorphism $\tau$ fixes
no vertices, no darts, and no directed and undirected edges, but some halvable edges may be fixed
(while transposing the corresponding pairs of darts).

\heading{Action of Automorphisms on Atoms.}
The following lemma explains how automorphisms behave with respect to the atoms of a graph.

\begin{lemma} \label{lem:atom_automorphisms}
Let $A$ be an atom and let $\pi \in \Aut(G)$. Then the following
holds:
\begin{packed_enum}
\item[(a)] The image $\pi(A)$ is an atom isomorphic to $A$. Further $\pi(\bo A) = \bo
\pi(A)$ and $\pi(\int A) = \int \pi(A)$, so $A \cong_\bo \pi(A)$.
\item[(b)] If $\pi(A) \ne A$, then $\pi(\int A) \cap \int A = \emptyset$.
\item[(c)] If $\pi(A) \ne A$, then $\pi(A) \cap A = \bo A \cap \bo \pi(A)$.
\end{packed_enum}
\end{lemma}

\begin{proof}
(a) Every automorphism permutes separately the set of articulations and non-trivial 2-cuts. So
$\pi(\bo A)$ separates $\pi(\int A)$ from the rest of the graph. It follows that $\pi(A)$ is an
atom, and $\pi$ clearly preserves the boundaries and the interiors.

For the rest, (b) follows from Lemma~\ref{lem:nonintersecting_interiors} and (c) follows from
Lemma~\ref{lem:nonintersecting_atoms}.\qed
\end{proof}

Therefore, every automorphism $\pi \in \Aut(G)$ induces a permutation of atoms, and $\Aut(G)$
induces an action on the set of all atoms.

\subsection{Regular Projections and Quotients of Atoms.} \label{sec:quotients_of_atoms}

Let $\Gamma$ be a semiregular subgroup of $\Aut(G)$, which defines a regular covering projection $p
: G \to G / \Gamma$. Negami~\cite[p.~166]{negami} investigated possible projections of proper atoms,
and we study this question in more detail. Let $A$ be a proper atom or a dipole with $\bo A =
\{u,v\}$ and $p |_A$ be its regular projection projection. We define three different types of projections $p
|_A$, illustrated in Fig.~\ref{fig:atom_projection_cases}:

\begin{figure}[b!]
\centering
\includegraphics{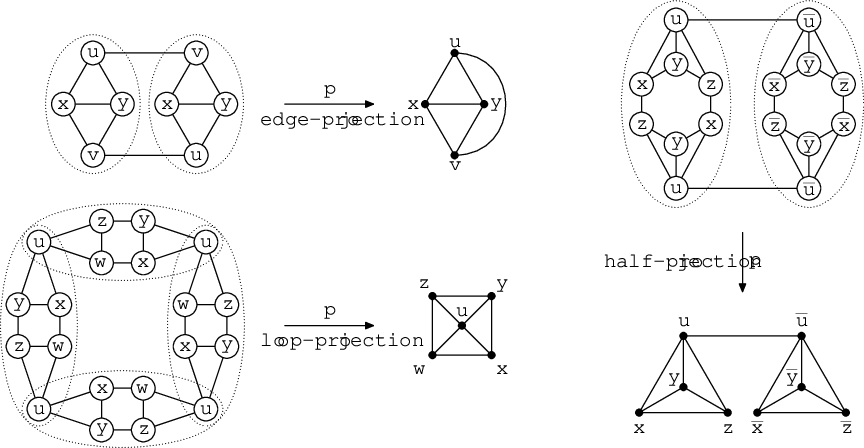}
\caption{Examples of three types for projections of atoms. Notice that for the third graph, an
	edge-projection can also be applied which gives a different quotient.}
	\label{fig:atom_projection_cases}
\end{figure}

\begin{packed_itemize}
\item \emph{An edge-projection $p |_A$.} The atom $A$ is preserved in $G / \Gamma$, meaning $p(A) \cong A$.
Notice that $p(A)$ is a subgraph of $G / \Gamma$, not necessarily induced. For instance for a proper
atom, it can happen that $p(u)p(v)$ is adjacent, even through $uv \notin \bE(G)$, as in
Fig.~\ref{fig:atom_projection_cases}.
\item \emph{A loop-projection $p |_A$.} The interior $\int A$ is preserved and the vertices $u$ and $v$ are
identified, i.e., $p(\int A) \cong \int A$ and $p(u) = p(v)$.
\item \emph{A half-projection $p |_A$.} The restriction of the covering projection
$p$ onto $A$ is a 2-fold covering with $\bo A$ forming a fibre. In particular, there exists an
involutory permutation $\pi$ in $\Gamma$ which exchanges $u$ and $v$ and preserves $A$. Then $p(u) =
p(v)$ and $p(A)$ consists of orbits of $\pi$ on $A$. This means that pairs of vertices and darts
of $A$ are identified in $p(A)$, so $\bv(p(A)) = \bv(A)/2$ and $\bd(p(A)) = \bd(A)/2$.
\end{packed_itemize}

\noindent For a block atom $A$, since $|\bo A| = 1$, we define only an \emph{edge-projection} $p
|_A$, when $p(A) \cong A$. 

\begin{lemma} \label{lem:atom_covering_cases}
Let $G$ be a graph, $A$ be its atom, and $p$ be a regular covering projection defined on $G$. Then
the restriction $p |_A$ is either an~edge-projection, or a~loop-projection, or it is
a~half-projection.  If $A$ is a block atom, then $p |_A$ is an edge-projection.  If $p |_A$ is a
half-projection, then $A$ is a halvable atom, and $p$ is a $2\ell$-fold regular covering projection
for some $\ell \in {\mathbb N}$.  
\end{lemma}

\begin{proof}
Let $\Gamma$ be the semiregular subgroup of $\Aut(G)$ defining $p : G \to G / \Gamma$, and we assume
that $\Gamma$ is non-trivial, otherwise $G \cong G / \Gamma$ and we get an edge-projection for all atoms.

We first deal with a block atom $A$. By Lemma~\ref{lem:central_block}, there exists a central block
$C$. In the notation of Lemma~\ref{lem:semiregular_action_on_blocks}, $A \subseteq G_u$ for some
articulation $u$ of $C$. By Lemma~\ref{lem:semiregular_action_on_blocks}, $G / \Gamma$ contains a
vertex $[u]$ together with a graph $p(G_u) \cong_\bo G_u$ attached. Therefore, $p|_A$ is an
edge-projection. 

It remains to deal with $A$ being a proper atom or a dipole, and let $\bo A = \{u,v\}$. According to
Lemma~\ref{lem:atom_automorphisms}b every automorphism $\pi$ either preserves $\int A$, or $\int A$
and $\pi(\int A)$ are disjoint. 

Suppose that there exists a non-trivial automorphism $\pi \in \Gamma$ preserving $\int A$. By
Lemma~\ref{lem:atom_automorphisms}a, we know $\pi(\bo A) = \bo A$, and by semiregularity, $\pi$ is
uniquely determined and exchanges $u$ and $v$.  Then the fiber containing $u$ and $v$ has to be of
an even size, with $\pi$ being an involution reflecting $\ell$ copies of $A$, and so $p$ is a $2\ell$-fold
covering projection. Therefore, $p |_A$ is a half-projection.

Suppose that there is no non-trivial automorphism which preserves $\int A$. The only difference
between an edge- and a loop-projection is whether $u$ and $v$ are contained in one fiber of $p |_A$,
or not. First, suppose that for every non-trivial $\pi \in \Gamma$ we get $A \cap \pi(A) =
\emptyset$.  Then no fiber contains more than one vertex of $A$, and $p |_A$ is an edge-projection,
i.e, $A \cong p(A)$. Next, suppose that there exists $\pi \in \Gamma$ such that $A \cap \pi(A) \ne
\emptyset$. By Lemma~\ref{lem:atom_automorphisms}c, we get $A \cap \pi(A) = \bo A \cap \bo \pi(A)$,
so $u$ and $v$ belong to one fiber of $p$, which makes $p |_A$ a loop-projection.\qed
\end{proof}

\begin{figure}[t!]
\centering
\includegraphics{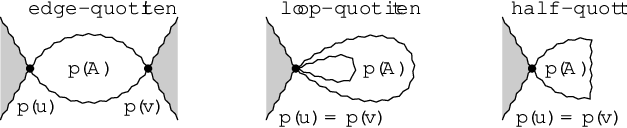}
\caption{How can the quotient $p(A)$ look in $G / \Gamma$, depending on type of $p |_A$.}
\label{fig:projections_of_atoms}
\end{figure}

Depending on the kind of the projection of an atom $A$ we call the respective  quotient $p(A)$:
\begin{packed_itemize}
\item an \emph{edge-quotient} if $p |_A$ is an edge-projection,
\item a \emph{loop-quotient} if $p |_A$ is a loop-projection,
\item a \emph{half-quotient} if $p |_A$ is a half-projection,
\end{packed_itemize}
\noindent see Figure~\ref{fig:projections_of_atoms} for an illustration. 

%For an edge-, a loop- and a half-projection $p |_A$, we get three types of quotients $p(A)$ of $A$
%which we call an \emph{edge-quotient}, a \emph{loop-quotient} and a \emph{half-quotient},
%respectively.

The following lemma allows to say ``the'' edge- and ``the'' loop-quotient of an atom since they are
uniquely determined.

\begin{lemma} \label{lem:unique_quotients}
Let $A$ be an atom and let $p$ and $p'$ be two regular covering projections. If $p |_A$ and $p' |_A$ are
both edge-projections, then $p(A) \cong_\bo p'(A)$. Similarly, if $A$ is a proper atom or a
dipole, and $p |_A$ and $p' |_A$ are both loop projections, then $p(A) \cong_\bo p'(A)$.
\end{lemma}

\begin{proof}
In both cases, we have $\int A \cong p(\int A)$ and $\int A \cong p'(\int A)$, so these quotients
are uniquely determined.\qed
\end{proof}

For half-quotients, this uniqueness does not hold. First, a proper atom or a dipole $A$ with $\bo A
= \{u,v\}$ has to be halvable to admit a half-quotient. Then each half-quotient $p(A)$ is determined
by a semiregular involutory automorphism $\tau \in \Autbo{A}$ exchanging $u$ and $v$. It is clear
that different automorphisms $\tau$ may give rise to non-isomorphic half-quotients. So when $p$ and
$p'$ are regular covering projections such that $p |_A$ and $p' |_A$ are half-projections, the
quotients $p(A)$ and $p'(A)$ might, or might not be $\bo$-isomorphic. In particular, for a dipole we
have the following sharp upper bound for the number of half-quotients.

\begin{lemma} \label{lem:dipole_quotients_noncolored}
A dipole $A$ has at most $\bigl\lfloor{\be(A) \over 2}\bigr\rfloor+1$ pairwise non-isomorphic
half-quotients, and this bound is achieved.
\end{lemma}

\begin{figure}[t!]
\centering
\includegraphics{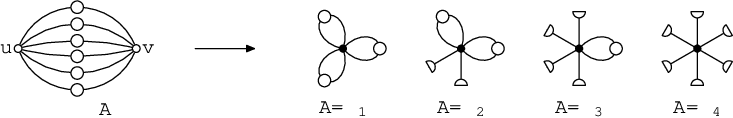}
\caption{Assuming that quotients can contain half-edges, the depicted dipole $A$ has four
non-isomorphic half-quotients $A / \Gamma_1,\dots,A / \Gamma_4$. We have $\Gamma_i = \langle \tau_i
\rangle$ where $\tau_i \in \Aut(\bo A)$ is an involution swapping $u$ and $v$, fixing exactly $2(i-1)$
edges and swapping pairs of the remaining edges.}
\label{fig:halfquotients_of_dipole}
\end{figure}

\begin{proof}
In order to prove the upper bound, we can assume  without loss of generality, that all edges of the
dipole $A$ are halvable. Let $\tau$ be a semiregular involution giving rise to the covering $p: A\to
A/\langle\tau\rangle$. The edges of $A$ that are fixed by $\tau$ project to half-edges in the
half-quotient $A / \left<\tau\right>$ while the pairs of edges of $A$ swapped by $\tau$ give rise to
loops in $A / \left<\tau\right>$. In the quotient, we have $\ell$ loops and $h$ half-edges attached
to a single vertex such that $2\ell+h = \be(A)$. Since $\ell$ ranges between $0$ and
$\bigl\lfloor{\be(A) \over 2}\bigr\rfloor$, the upper bound is established.

Finally, for each pair $(h,\ell)$ of integers satisfying the equation $2\ell+h = \be(A)$ one can
find an involution $\tau_h$  acting on the $2\be(A)$ darts of $A$, swapping the two vertices such
that the induced involution on the edges has $\ell$ orbits of size two and $h$ fixed points.
Employing $\tau_h$ for all admissible $h$, we get the  $\bigl\lfloor{\be(A) \over 2}\bigr\rfloor+1$
pairwise non-isomorphic quotients $A/\langle\tau_h\rangle$.
Figure~\ref{fig:halfquotients_of_dipole} shows the construction for $\be(A)=6$. \qed
\end{proof}

In contrast, for planar proper atoms, we prove in Lemma~\ref{lem:planar_proper_half_quotients} that there are at
most two non-isomorphic half-quotients. The non-uniqueness of half-quotients is one of the main
algorithmic difficulties for an effective regular covering testing in the class of planar graphs in~\cite{fkkn,fkkn16}.

\section{Graph Reductions and Quotient Expansions} \label{sec:reduction_and_expansion}

We start with a quick overview of the reduction procedure. The reduction initiates with a graph $G$
and produces a sequence of graphs $G = G_0,G_1,\dots,G_r$. To produce $G_{i+1}$ from $G_i$, we find
the collection of all atoms $\calA$ in $G_i$ and replace each of them by an edge of the
corresponding type and color. We stop after $r$ steps when a primitive graph $G_r$ containing no
further atoms is reached. We call this sequence of graphs starting with $G$ and ending with a
primitive graph $G_r$ the \emph{reduction series} of $G$.

In this section, we describe structural properties of the reductions and expansions. We study
changes of the automorphism groups in each step of reduction series. The reduction forming $G_{i+1}$
from $G_i$ is defined in a way that an essential information of $\Aut(G_i)$ is preserved in
$\Aut(G_{i+1})$. In particular, we show that a semiregular subgroup $\Gamma=\Gamma_0$ gives rise to
a uniquelly determined sequence $\Gamma_0,\Gamma_1,\ldots,\Gamma_r$ of semiregular groups of
automorphisms such that $\Gamma_i\leq \Aut(G_i)$ and $\Gamma_i\cong \Gamma$, for $i=0,1,\ldots,r$.
This is a key property of the reductions which allows to realize the backward process of the
expansions starting from the covering $G_r\to G_r/\Gamma_r$, in case the regular covering $G\to
H=G/\Gamma$ exists.  More precisely, suppose that $H_r = G_r / \Gamma_r$ is a quotient of $G_r$ for
some semiregular group $\Gamma_r$. The reductions applied to reach $G_r$ are reverted on $H_r$ and
produce an \emph{expansion series} $H_r,H_{r-1},\dots,H_0$ of $H_r$.  Processing the expansion we
obtain a series of semiregular subgroups $\Gamma_r,\dots,\Gamma_0$ such that $H_i = G_i / \Gamma_i$
and $\Gamma_i$ extends the action of $\Gamma_{i+1}$. Thus each regular covering defined on $G$ can
be reconstructed from a regular covering defined on the associated primitive graph $G_r$. The entire
process is depicted in the diagram in Fig.~\ref{fig:reduction_expansion_diagram}.

The problem is that the expansions are, unlike the reductions, not uniquely determined.
From $H_{i+1}$, we can construct multiple $H_i$. In this section, we characterize all possible ways
how $H_i$ can be constructed from $H_{i+1}$ thus establishing Theorem~\ref{thm:quotient_expansion}.

\subsection{Reducing Graphs Using Atoms} \label{sec:reduction}

The reduction produces a series of graphs $G = G_0,\dots,G_r$, by replacing atoms with colored edges
encoding isomorphism classes. The edges are endowed with the edge types encoding the symmetry types of atoms.

\begin{quote}
{\bf Remark:}\enspace\textit{In what follows, we work with multigraphs with colored edges of the three
types: halvable, undirected and directed. For every automorphism/isomorphism, we require that it
preserves colors, edge types and direction of oriented edges.}
\end{quote}

\noindent We note that the results established in Section~\ref{sec:atoms} extend to colored graphs and
colored atoms without any problems.

\begin{figure}[t!]
\centering
\includegraphics{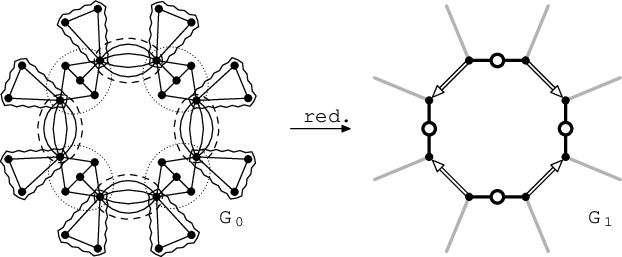}
\caption{On the left, we have a graph $G_0$ with three isomorphism classes of atoms, each having
four atoms.  The dipoles are halvable, the block atoms are symmetric and the proper atoms are
asymmetric.  We reduce $G_0$ to $G_1$ which is a cycle with attached single pendant edges, with four
black halvable edges replacing the dipoles, eight gray pendant edges replacing the block atoms,
and four white directed edges replacing the proper atoms. The reduction series ends with $G_1$ since
it is primitive. Notice the consistent orientation of the directed edges.}
\label{fig:example_of_reduction}
\end{figure}

For a graph $G_i$, we find the collection of all atoms $\calA$. We obtain $\bo$-isomorphism
classes for the set of all atoms $\calA$ of $G_i$ such that $A$ and $A'$ belong to the same class if
and only if $A \cong_\bo A'$. To each isomorphism class, we assign one new color not yet used in
the process of the reduction. The graph $G_{i+1}$ is constructed from $G_i$ by replacing the
interior of each atom in $\calA$ by an edge as follows:
\begin{packed_itemize}
\item \emph{The interior of a block atom $A$} with $\bo A = \{u\}$ is replaced by a pendant edge
based at $u$ of the color assigned to the isomorphism class containing $A$.
\item \emph{The interior of a proper atom or a dipole $A$} with $\bo A = \{u,v\}$, which is
halvable/symmetric/asymmetric, is replaced by a new halvable/undirected/directed edge $uv$,
respectively, of the color assigned to the isomorphism class containing $A$.  The orientation of the
edges is defined as follows.  For each isomorphism class of asymmetric atoms, we choose an arbitrary
orientation of the directed edge replacing an atom $A$ in the class. If $A'$ is another atom in the
class, then each $\bo$-isomorphism takes $\bo A\mapsto \bo A'$ uniformly thus prescribing the
orientation of the edge replacing $A'$ in the reduction. In particular, the edges replacing
asymmetric atoms in the same orbit of $\Aut(G_i)$ are oriented consistently with the action of the
group.
\end{packed_itemize}
For example of the reduction, see Fig.~\ref{fig:example_of_reduction}.

\begin{wrapfigure}{r}{0.45\textwidth}
\begin{center}
\vskip -2ex
\includegraphics[width=0.45\textwidth]{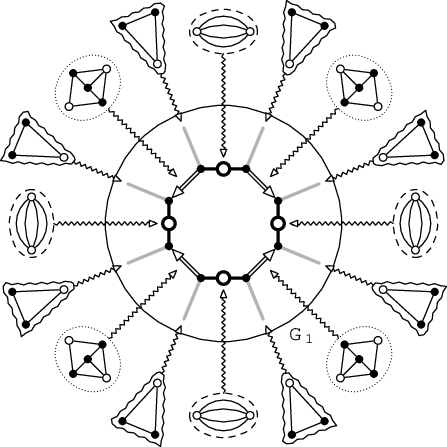}
\end{center}
\vskip -5ex
\caption{The reduction tree for the reduction series in Fig.~\ref{fig:example_of_reduction}. The
root is the primitive graph $G_1$ and each leaf corresponds to one atom of $G_0$. For all atoms $A$, the
vertices of $\bo A$ are depicted in white.}
\label{fig:reduction_tree}
\vskip -5ex
\end{wrapfigure}

According to Lemma~\ref{lem:nonintersecting_atoms}, the replaced interiors of the atoms of $\calA$
are pairwise disjoint, so the reduction is well defined. We stop in the step $r$ when $G_r$ is a
primitive graph containing no atoms. (Recall Lemma~\ref{lem:primitive_graphs} characterizing all
primitive graphs.)

\begin{figure}[p!]
\centering
\includegraphics[width=\textwidth]{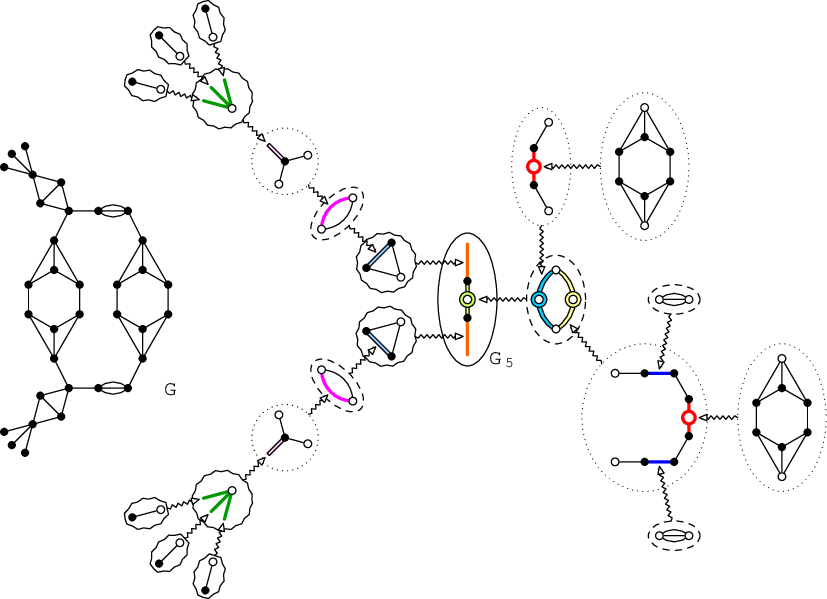}
\caption{A graph $G$ on the left and its reduction tree on the right, having the primitive graph
$G_5$ as the root. For all atoms $A$, the vertices of $\bo A$ are depicted in white.  The reader may
try to deduce the reduction series $G = G_0,G_1,\dots,G_5$ from this reduction tree.}
\label{fig:reduction_tree_large}
\end{figure}

For every graph $G$, the reduction series determines  the \emph{reduction tree} $T$ which is a
rooted tree defined as follows.  The root is the primitive graph $G_r$, and the other nodes are the
atoms obtained during the reductions. If a node $N$ of $T$ contains a colored edge, the
corresponding atom $A$ is a child of $N$ in $T$. Therefore, the leaves of $T$ are the atoms of
$G_0$, after removing them, the new leaves are the atoms of $G_1$, and so on. For examples, see
Fig.~\ref{fig:reduction_tree} and~\ref{fig:reduction_tree_large}.

\begin{figure}[p!]
\centering
\includegraphics{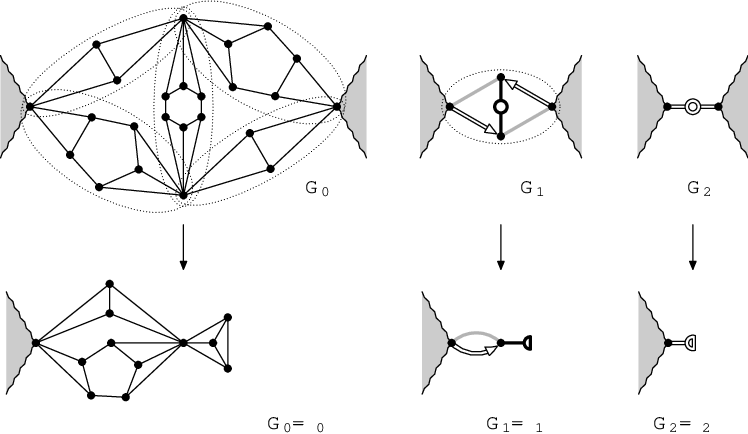}
\caption{We reduce a part of a graph in two steps. In the first step, we replace five atoms by five
edges. As the result we obtain one halvable atom which we further reduce to one halvable edge.
Notice that without considering edge types, the resulting atom in $G_1$ would be just symmetric. In
the bottom, we show a part of the corresponding quotient graphs when $\Gamma_i$ contains a
semiregular involutory automorphism $\pi$ from the definition of a half-projection.}
\label{fig:reducing_atoms}
\end{figure}

Fig.~\ref{fig:reducing_atoms}  shows a regular 2-fold covering $G_0\to G_0/\Gamma_0$, where
$\Gamma_0$ is a semiregular group of order two rotating the central hexagon. Observe that although
none of the graphs $G_0$ and $H_0=G_0/\Gamma_0$ contain half-edges,  in the quotients $G_1 /
\Gamma_1$ and $G_2 / \Gamma_2$  half-edges do appear. This example shows that in the reductions and
expansions we need to consider half-edges even if both graphs $G$ and $G / \Gamma$ are simple.

\subsection{Reduction Epimorphism} \label{sec:reduction_epimorphism}

We describe algebraic properties of the reductions, in particular how the groups $\Aut(G_i)$ and
$\Aut(G_{i+1})$ are related. There exists a natural mapping $\Phi_i: \Aut(G_i) \to \Aut(G_{i+1})$
called \emph{reduction epimorphism} which we define as follows. Let $\pi \in \Aut(G_i)$. For the
common vertices, darts and edges of $G_i$ and $G_{i+1}$, we define the action of $\Phi_i(\pi)$
following the action of  $\pi$.  If $A$ is an atom of $G_i$, then according to
Lemma~\ref{lem:atom_automorphisms}a, $\pi(A)$ is an atom such that $\pi(A) \cong_\bo A$. In
$G_{i+1}$, we replace, respectively, the interiors of both $A$ and $\pi(A)$ by the edges $e_A$ and
$e_{\pi(A)}$ of the same type and color. We define $\Phi_i(\pi)(e_A) = e_{\pi(A)}$.  By
Lemma~\ref{lem:atom_automorphisms},  $\Phi_i(\pi)$ is a well-defined automorphism
of $G_{i+1}$.

For purpose of Section~\ref{sec:quotients}, we also define $\Phi_i$ on the darts.  We choose $u \in
\bo A$ and $\pi(u) \in \bo \pi(A)$. Let $d_u$ be the dart of $e_A$ incident with $u$ and $d'_u$ be
the other dart, and similarly let $d_{\pi(u)}$ and $d'_{\pi(u)}$ be the darts composing
$e_{\pi(A)}$.  Then we define $\Phi_i(\pi)(d_u) = d_{\pi(u)}$ and $\Phi_i(\pi)(d'_u) = d'_{\pi(u)}$.

\begin{proposition} \label{prop:reduction_homomorphism}
The mapping $\Phi_i : \Aut(G_i) \to \Aut(G_{i+1})$, $0\leq i<r$, satisfies the following:
\begin{packed_enum}
\item[(a)] The mapping $\Phi_i$ is a group homomorphism.
\item[(b)] The mapping $\Phi_i$ is an epimorphism, i.e., it is surjective.
\item[(c)] For a semiregular subgroup $\Gamma$ of $\Aut(G_i)$, the restriction $\Phi_i |_\Gamma$ is an
isomorphism. Moreover, the subgroup $\Phi_i(\Gamma)$ remains semiregular in $\Aut(G_{i+1})$.
\end{packed_enum}
\end{proposition}

\begin{proof}
(a) Clearly, $\Phi_i(\id) = \id$. Let $\pi,\sigma \in \Aut(G_i)$. We need to show that
$\Phi_i(\sigma \pi) = \Phi_i(\sigma) \Phi_i(\pi)$. This is clearly true outside the interiors of the
atoms. Let $A$ be an atom. By the definition, $\Phi_i(\sigma \pi)$ maps $e_A$ to
$e_{\sigma(\pi(A))}$ while $\Phi_i(\pi)$ maps $e_A$ to $e_{\pi(A)}$ and $\Phi(\sigma)$ maps
$e_{\pi(A)}$ to $e_{\sigma(\pi(A))}$. So the equality holds everywhere and $\Phi_i$ is a group
homomorphism. 

(b) Let $\pi' \in \Aut(G_{i+1})$, we want to extend $\pi'$ to $\pi \in \Aut(G_i)$ such that
$\Phi_i(\pi) = \pi'$. We just describe this extension on a single edge $e$. If $e$ is
an original edge of $G_i$, there is nothing to extend. Suppose that $e$  in $G_{i+1}$ replaced
an atom $A$ in $G_i$. Then $\hat e = \pi'(e)$ is an edge of the same color and the same type as $e$, and
therefore $\hat e$  comes from an  atom $\hat A$ isomorphic to $A$ of the same symmetry type. The
automorphism $\pi'$ prescribes the action on the boundary $\bo A$. We need to show that it is
possible to extend the action on $\int A$ consistently. We distinguish three cases.
\begin{packed_itemize}
\item \emph{$A$ is a block atom:} The edges $e$ and $\hat e$ are pendant, attached by articulations
$u$ and $u'$. We define $\pi |_A$ equal to an arbitrary $\bo$-isomorphism from $A$ to $\hat
A$.
\item \emph{$A$ is an asymmetric proper atom/dipole:} Suppose that $e = uv$ and $\hat e = \hat u
\hat v$ are oriented from $u$ to $v$ and from $\hat u$ to $\hat v$. We define $\pi |_A$ equal to an
arbitrary $\bo$-isomorphism $\sigma$ from $A$ to $\hat A$. Since the orientation of $e$ and $\hat e$
is consistent, we have that $\sigma(u) = \hat u$ and $\sigma(v) = \hat v$, so the mappings $\pi$ and
$\pi'$ are the same on $\bo A$.
\item \emph{$A$ is a symmetric/halvable proper atom/dipole:} Let $\sigma$ be an arbitrary
$\bo$-isomorphism from $A$ and $\hat A$. Either $\sigma$ maps $\bo A$ exactly as $\pi'$, and then we
define $\pi |_A = \sigma$. Otherwise let $\tau \in \Autbo{A}$ be an automorphism exchanging the two
vertices of $\bo A$ which exists since $A$ is not asymmetric. Then we define $\pi |_A = \tau
\sigma$.
\end{packed_itemize}
So $\Phi_i$ is a surjective mapping.

(c) Recall that the kernel $\Ker(\Phi_i)$ is the set of all $\pi$ such that $\Phi_i(\pi) = \id$ and
it is a normal subgroup of $\Aut(G_i)$. It has the following structure: $\pi \in \Ker(\Phi_i)$ if
and only if it fixes everything except for the interiors of the atoms. Further, $\pi(\int A) = \int
\pi(A)$, so $\pi\in\Ker(\Phi_i)$ can non-trivially act only inside the interiors of the atoms. 

For any subgroup $\Gamma$, the restricted mapping $\Phi_i |_\Gamma$ is a group homomorphism with
$\Ker(\Phi_i |_\Gamma) = \Ker(\Phi_i) \cap \Gamma$. If $\Gamma$ is semiregular, then we show that
$\Ker(\Phi_i) \cap \Gamma$ is trivial. We know that $G_i$ contains at least one atom $A$. The
boundary $\bo A$ is point-wise fixed by $\Ker(\Phi_i)$, so by the semiregularity of $\Gamma$ the
intersection $\Gamma\cap\Ker(\Phi_i)$ is trivial. Hence $\Phi_i |_\Gamma$ is an isomorphism.

For the semiregularity of $\Phi_i(\Gamma)$, let $\pi' \in \Phi_i(\Gamma)$. Since $\Phi_i|_\Gamma$ is
an isomorphism, there exists the unique $\pi \in \Gamma$ such that $\Phi_i(\pi) = \pi'$.  Recall
that the definition of semiregularity requires trivial stabilizers of vertices, darts, and
non-halvable edges, but admits non-trivial stabilizers of halvable edges. If $\pi'$ fixes a vertex
$u$, then $\pi$ fixes $u$ as well, so necessarily $\pi = \id$ and $\pi' = \Phi_i(\id) = \id$. A
non-trivial semiregular automorphism $\pi'$ might fix an edge $e = uv$ while exchanging $u$ and $v$.
We need to verify that $e$ is halvable.  If $e$ belongs to both $G_i$ and $G_{i+1}$, then $\pi$ also
fixes $e$, so $e$ is halvable and $\pi'$ can fix it as well.  Otherwise there is an atom $A$ in
$G_i$ replaced by $e$ in $G_{i+1}$. Then $\pi|_A$ is an involutory semiregular automorphism
exchanging $u$ and $v$, so $A$ is halvable. But then $e$ is a halvable edge, and thus $\pi'$ is
allowed to fix it.\qed
\end{proof}

The above statement is an example of a phenomenon known in permutation group theory. Interiors
of atoms behave as \emph{blocks of imprimitivity} in the action of $\Aut(G_i)$. It is well-known
that the kernel of the action on the imprimitivity blocks is a normal subgroup of $\Aut(G_i)$.

Now, we are ready to prove Proposition~\ref{prop:group_extension} which states that $\Aut(G_{i+1})
\cong \Aut(G_i) / \Ker(\Phi_i)$:

\begin{proof}[Proposition~\ref{prop:group_extension}]
By Proposition~\ref{prop:reduction_homomorphism}b, $\Phi_i$ is surjective, so by the well-known
1st Isomorphism Theorem it follows that $\Aut(G_{i+1}) \cong \Aut(G_i) / \Ker(\Phi_i)$.\qed
\end{proof}

In~\cite{knz}, the opposite relation between $\Aut(G_i)$ and $\Aut(G_{i+1})$ is studied. Under the
assumption that every symmetric proper atom has an involutory automorphisms $\tau \in \Autbo{A}$
exchanging the vertices of $\bo A$, the argument in the proof of
Proposition~\ref{prop:reduction_homomorphism}b is expanded and it is shown that $\Aut(G_i) \cong
\Aut(G_{i+1}) \ltimes \Ker(\Phi_i)$.

\begin{corollary}
We have $\Aut(G_r) \cong \Aut(G_0) / \Ker(\Phi_{r-1} \circ \Phi_{r-2} \circ \cdots \circ \Phi_0)$.
\end{corollary}

\begin{proof}
We prove by induction that $\Aut(G_i) \cong \Aut(G_0) / \Ker(\Phi_{i-1} \circ \Phi_{i-2} \circ
\cdots \circ \Phi_0)$ where the first step holds by Proposition~\ref{prop:group_extension}. For the
induction step, we use the well-known 3rd Isomorphism Theorem stating that for normal subgroups
$\Psi \le \Psi'$ of $\Sigma$, it holds that $(\Sigma / \Psi) / (\Psi' / \Psi) \cong \Sigma / \Psi'$.
We choose
$$\Sigma = \Aut(G_0),\quad \Psi = \Ker(\Phi_{i-2} \circ \cdots \Phi_0),\quad \text{and}\quad
\Psi' = \Ker(\Phi_{i-1} \circ \cdots \circ \Phi_0).$$
By the induction hypothesis, $\Sigma / \Psi \cong \Aut(G_{i-1})$. Since $\Psi' / \Psi \cong
\Ker(\Phi_{i-1})$, we get that $\Aut(G_{i-1}) / \Ker(\Phi_{i-1}) \cong \Aut(G_i) \cong \Aut(G_0) /
\Ker(\Phi_{i-1} \circ \cdots \circ \Phi_0)$ (where the first part holds by
Proposition~\ref{prop:group_extension}).
\qed
\end{proof}

We can also describe the structure of $\Ker(\Phi_i)$:

\begin{lemma} \label{lem:kernel_char}
We have
$$\Ker(\Phi_i) \cong \prod_{A \in \calA} \Fix(\bo A).$$
\end{lemma}

\begin{proof}
Every automorphism of $\Ker(\Phi_i)$ fixes all vertices and edges outside the interiors of atoms of
$\calA$. According to Lemma~\ref{lem:nonintersecting_interiors}, these interiors are pairwise
disjoint, so $\Ker(\Phi_i)$ acts independently on each interior. Thus we get $\Ker(\Phi_i)$ as the
direct product of actions on each interior $\int A$ which is precisely $\Fix(\bo A)$.\qed
\end{proof}

Let $A_1,\dots,A_s$ be pairwise non-isomorphic atoms in $G_i$, appearing with multiplicities
$m_1,\dots,m_s$. According to Lemma~\ref{lem:kernel_char}, we get
$$\Ker(\Phi_i) \cong \Fix(A_1)^{m_1} \times \cdots \times \Fix(A_s)^{m_s}.$$
For the example of Fig.~\ref{fig:example_of_reduction}, we have $\Ker(\Phi_0) \cong \gC_2^8 \times
\gC_2^4 \times \gS_4^4$. For the example in Fig.~\ref{fig:example_of_reduction}, it is shown
in~\cite[Proposition 4.7]{knz} that
$$\Aut(G_1) \cong \gC_2^2\qquad\text{and}\qquad
\Aut(G_0) \cong (\gC_2^8 \times \gC_2^4 \times \gS_4^4) \rtimes \gC_2^2.$$

\subsection{Reduction Preserves the Central Block}

We show that the reduction preserves the central block:

\begin{lemma} \label{lem:preserved_center}
Let $G$ admit a non-trivial semiregular automorphism $\pi$. Then each $G_{i+1}$ has a central block
which is obtained from the central block of $G_i$ by replacing its proper atoms and dipoles by colored edges.
\end{lemma}

\begin{proof}
By Proposition~\ref{prop:reduction_homomorphism}c, semiregular automorphisms are preserved during
the reduction. Let $\pi_i$ be a non-trivial semiregular automorphism of $G_i$, such that $\pi_0 =
\pi$ and $\pi_{i+1} = \Phi_i(\pi_i)$. By Lemma~\ref{lem:central_block}, each $G_i$ has a central
block. Since we replace only proper atoms and dipoles in the central block, it remains to be a block
after the reduction. We argue by induction that it remains central as well.

Let $C$ be the central block of $G_i$ and let $C'$ be the block in $G_{i+1}$ obtained by reducing
all atoms of $C$. Let $T$ and $T'$ be the block trees of $G_i$ and $G_{i+1}$, respectively, and we
assume that both are rooted towards $C$ and $C'$, respectively (even though we still need to prove
that $C'$ is the central block of $G_{i+1}$).  Let $u$ be an articulation of $C'$ such that the
rooted subtree $T'_u$ of $T'$ determined by $u$ contains a longest oriented path towards $u$. Since
$u$ also belongs to $C$, let $T_u$ be the corresponding rooted subtree of $T$ determined by $u$. 

Since $\pi_i$ is semiregular, we have $\pi_i(u) = v$ for some $v \ne u$. The vertex $v$ is an
articulation of $G_i$ and let $T_v$ be a rooted subtree of $T$ determined by $v$. By
Lemma~\ref{lem:semiregular_action_on_blocks}, we have $G_v \cong G_u$ where $G_u$ and $G_v$ are the
subgraphs of $G_i$ induced by $T_u$ and $T_v$, respectively.  In $G_{i+1}$, we apply one step of
reductions on both $G_u$ and $G_v$, and we obtain $G'_u \cong G'_v$, both attached to $u$ and $v$ in
$C'$, respectively. Therefore, $C'$ is the central block of $G_{i+1}$.\qed
\end{proof}

Recall that the definition of atoms depends on the central block.  The above statement establishes
that for graphs admitting a non-trivial semiregular automorphism ``the position'' of the central
block does not change when applying the reductions. Therefore, the primitive graph $G_r$ contains a
central block, and by Lemma~\ref{lem:primitive_graphs} it is either 3-connected, or $C_n$, or $K_2$,
or can be made from these graph by attaching single pendant edges to at least two vertices.  In
general, the central block/articulation does not have to be preserved in the reduction and one has
to define atoms in all steps of the reduction with respect to the ``images'' of the same original
block/articulation. 

\subsection{Quotients and Their Expansion} \label{sec:quotients}

Let $G_0,\dots,G_r$ be the reduction series of $G$ and let $\Gamma_0$ be a semiregular subgroup
of $\Aut(G_0)$. By repeated application of Proposition~\ref{prop:reduction_homomorphism}c, we get
the uniquely determined semiregular subgroups $\Gamma_1,\dots,\Gamma_r$ of
$\Aut(G_1),\dots,\Aut(G_r)$ such that $\Gamma_{i+1} = \Phi_i(\Gamma_i)$, each isomorphic to
$\Gamma_0$. Let $H_i = G_i / \Gamma_i$ be the quotients with preserved colors of edges, and let
$p_i$ be the corresponding covering projection from $G_i$ to $H_i$.  Recall that $H_i$ can contain
edges, loops and half-edges. We summarize it as follows.

\begin{lemma} \label{lem:group_reduction}
Every semiregular subgroup $\Gamma_i$ of $\Aut(G_i)$ corresponds to a unique semiregular subgroup
$\Gamma_{i+1}$ of $\Aut(G_{i+1})$ such that $\Gamma_{i+1} = \Phi_i(\Gamma_i)$. Moreover,
$\Gamma_i \cong \Gamma_{i+1}$. \qed
\end{lemma}

\begin{figure}[b!]
\centering
\includegraphics{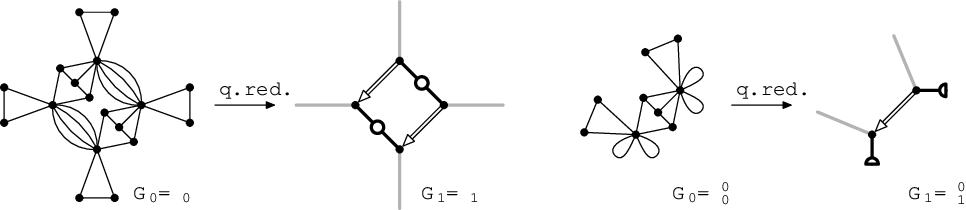}
\caption{An example of two quotients of the graph $G_0$ from Fig.~\ref{fig:example_of_reduction} with
the corresponding quotients of the reduced graph $G_1$. Here $\Gamma_1 = \Phi_1(\Gamma_0)$ and
$\Gamma'_1 = \Phi_1(\Gamma'_0)$. Both $\Gamma_0$ and $\Gamma_1$ are generated by the $180^\circ$
rotation of $G_0$ and $G_1$, respectively, while $\Gamma'_0$ and $\Gamma'_1$ are generated by
vertical and horizontal reflections of $G_0$ and $G_1$, where the reflections of $G_0$ swap pairs of
edges for dipoles whose boundaries exchange.}
\label{fig:computed_quotients}
\end{figure}

\heading{Quotient Reductions.}
Let $H_i = G_i / \Gamma_i$ and $p_i: G_i\to H_i$ be the regular covering induced by the action of a
semiregular group $\Gamma_i$. We investigate the relation between $H_i$ and $H_{i+1}$. Let $A$ be an
atom of $G_i$ represented by a colored edge $e$ in $G_{i+1}$.  According to
Lemma~\ref{lem:atom_covering_cases}, $p_i |_A$ is either an edge-, a loop-, or a half-projection. It
is easy to see that $\Gamma_{i+1}=\Phi_i(\Gamma_i)$ is defined exactly in the way that $p_{i+1}(e)$
is an edge if $p_i |_A$ is an edge-projection,  a loop if $p_i |_A$ is a loop-projection, and is a
half-edge if $p_i |_A$ is a half-projection. (This explains the choice of names for projections and
quotients of atoms.) Figure~\ref{fig:computed_quotients} shows examples.

We define the \emph{quotient reduction} of $H_i$ by replacing the projections of atoms $A$ of $G_i$
with the projections of the corresponding colored edges of $G_{i+1}$. More precisely, the
edge-quotient $p_i(A)$ in $H_i$ of an atom $A$ in $G_i$ is replaced by an edge of the same color and
type as the edge $e$ replacing $A$ in $G_i$, and similarly loop-quotients and half-quotients are
replaced by colored loops and half-edges. The quotient reduction produces from $H_i$ the graph
$H_{i+1}$ which is a quotient of $G_{i+1}$.  In other words, the quotient reduction is defined in
such a way that the following diagram commutes: 
\begin{equation} \label{eq:red_diagram}
\begin{gathered}
\includegraphics{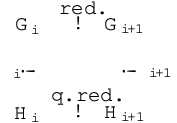}
\end{gathered}
\end{equation}

\heading{Overview of Quotient Expansions.}
Our goal is to reverse the horizontal edges in Diagram~(\ref{eq:red_diagram}), i.e, to understand:
\begin{equation} \label{eq:exp_diagram}
\begin{gathered}
\includegraphics{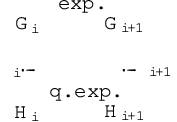}
\end{gathered}
\end{equation}
Here, the \emph{expansion} of $G_{i+1}$ replaces colored edges corresponding to atoms of $G_i$ by these atoms,
so it constructs $G_{i+1}$. As we will see, the bottom arrow is the \emph{quotient expansion} of
$H_{i+1}$ which replaces colored edges, loops and half-edges corresponding to quotients of atoms of
$G_i$ by edge-, loop-, and half-quotients of these atoms, so it constructs $H_i$.

Let $\Gamma_i$ and $\Gamma_{i+1}=\Phi_i(\Gamma_i)$ be semiregular groups of $\Aut(G_i)$ and $\Aut(G_{i+1})$,
respectively.    We call such $\Gamma_{i+1}$ a \emph{reduction} of $\Gamma_i$, and $\Gamma_i$ an
\emph{extension} of $\Gamma_{i+1}$. Now, we study the possible extensions of $\Gamma_{i+1}$, i.e,
all groups $\Gamma_i \le \Aut(G_i)$ such that $\Phi_i(\Gamma_i) = \Gamma_{i+1}$.

There are two fundamental questions we address in this section in full detail:
\begin{packed_itemize}
\item \emph{Question 1.} Given a group $\Gamma_{i+1}$, which semiregular groups $\Gamma_i$ are its
extensions? Notice that all these extensions of $\Gamma_{i+1}$ are isomorphic to $\Gamma_{i+1}$ as abstract
groups, but they may act differently on $G_i$.
\item \emph{Question 2.} Let $\Gamma_i$ and $\Gamma'_i$ be two semiregular groups extending $\Gamma_{i+1}$.
Under which conditions are the quotients $H_i = G_i / \Gamma_i$ and $H'_i = G_i / \Gamma'_i$
different?
\end{packed_itemize}

\heading{Extensions of Group Actions.} We first deal with Question 1.

\begin{lemma} \label{lem:group_expansion}
For every semiregular group $\Gamma_{i+1} \le \Aut(G_{i+1})$, there exists an extension $\Gamma_i
\le \Aut(G_i)$ such that $\Phi_i(\Gamma_i) = \Gamma_{i+1}$.
\end{lemma}

\begin{proof}
First notice that $\Gamma_{i+1}$ determines the action of $\Gamma_i$ everywhere on $G_i$ except for the
interiors of the atoms of $G_i$, so we just need to define it there. Let $e$ be one edge of
$G_{i+1}$ replacing an atom $A$ of $G_i$. Let $|\Gamma_{i+1}| = k$. Semiregularity of $\Gamma_{i+1}$
implies that the orbit $[e]$ is either of size $k$, or of size $k \over 2$. Let $\pi' \in
\Gamma_{i+1}$. To define its extension $\pi$ on the interiors of the atoms of $G_i$ we distinguish
three cases, see Fig.~\ref{fig:group_extension}:

\begin{figure}[t!]
\centering
\includegraphics{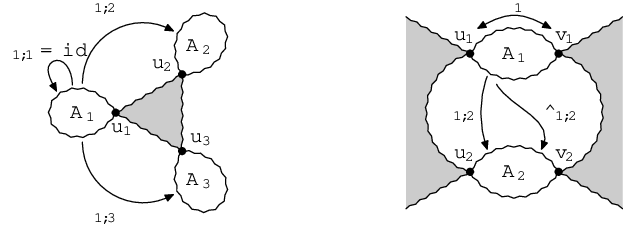}
\caption{Case 1 is depicted on the left for three edges corresponding to isomorphic block atoms $A_1$,
$A_2$ and $A_3$. The depicted isomorphisms are used to extend $\Gamma_{i+1}$ on the interiors of these
atoms. Case 3 is on the right, with an additional semiregular involution $\tau_1$ which transposes
$u_1$ and $v_1$.}
\label{fig:group_extension}
\end{figure}

\emph{Case 1: The orbit $[e]$ is of size $k$ and $e$ is a pendant edge.}
Then $A$ is a block atom in $G_i$.  Suppose that $e$ is attached to $u$.  Let $[e] =
\{e_1,\dots,e_k\}$, $[u]=\{u_1,u_2,\dots,u_k\}$. Denote by $\pi'_j$ the unique automorphism in
$\Gamma_{m+1}$ such that $u_j = \pi_j'(u_1)$ and  $\pi_j'(e)=e_j$. Note that
$\Gamma_{m+1}=\{\pi'_1,\pi'_2,\dots,\pi'_k\}$.

Let $A_1,\dots,A_k$ be the atoms of $G_i$ corresponding to the edges $e_1,\dots,e_k$ in $G_{i+1}$.
The edges $e_1,\dots,e_k$ have the same color and type, and thus the block atoms $A_j$, for
$j=1,2,\dots,k$, are pairwise isomorphic.

We define the action of $\Gamma_i$ on the interiors of $A_1,\dots,A_k$ as follows.
We choose arbitrarily $\bo$-isomorphisms $\sigma_{1,x}$ from $A_1$ to $A_x$, for $x = 2,\dots,k$, such that
$\sigma_{1,x}(u_1) = u_x$, and put $\sigma_{1,1} = \id$ and $\sigma_{x,y} =
\sigma_{1,y}\sigma^{-1}_{1,x}$. If $\pi'(e_x) = e_y$, we set $\pi |_{\int A_x} =
\sigma_{x,y} |_{\int A_x}$. Since
\begin{equation} \label{eq:sigma}
\sigma_{x,z} = \sigma_{y,z} \sigma_{x,y},\qquad \forall x,y,z \in \{1,\dots,k\},
\end{equation}
the composition of the extensions $\pi_x$ and $\pi_y$ of $\pi'_x$ and $\pi'_y$, respectively, is
defined on the interiors of $A_j$ for $j=1,\dots,k$ consistently. Hence
$\Gamma_i=\{\pi_1,\pi_2,\dots,\pi_k\}$ is a group of automorphisms acting on the partial extension,
where only the edges $e_j$ are expanded to the atoms $A_j$, for $j=1,2,\dots,k$. 

We show that the action remains semiregular. Suppose $\pi_j(v)=v$ for some vertex $v$. Then either
$v$ is a vertex of $G_{i+1}$, or it is an internal vertex of an atom $A_x$. The second case reduces
to the first one since the articulation $u_x$ is fixed by $\pi_j$. We see that $\pi_j$ is an
extension of the identity, so it maps $\int A_y$ to $\int A_y$, for $y=1,\dots,k$. Since $\pi_j
|_{\int A_y} = \sigma_{y,y} |_{\int A_y} = \id$, each vertex and dart of $A_i$ is fixed, so $\pi_j
= \id$. Next, suppose that non-trivial $\pi_j \in \Gamma_i$ fixes an edge $\hat e$, and we want to show that
$\hat e$ is halvable. If $\hat e$ is also an edge of $G_{i+1}$, it is fixed by non-trivial $\pi'_j$
as well, so it is halvable. Otherwise, $\hat e$ is an internal edge of $A_x$, and as argued above,
both $\pi_j$ and $\pi'_j$ are identities, contradicting non-triviality. Hence the extended action
$\Gamma_i$ remains semiregular.

\emph{Case 2: The orbit $[e]$ is of size $k$ and $e$ is not a pendant edge.}
Let $e = uv$, and it corresponds to a proper atom or dipole $A$ in $G_i$.  Let $[e] =
\{e_1,\dots,e_k\}$, where $e_j=u_jv_j$, and let $[u] = \{u_1,\dots,u_k\}$ and $[v] =
\{v_1,v_2,\dots,v_k\}$.  We denote by $\pi'_j$ the semiregular automorphism of $\Gamma_{i+1}$ such
that $u_j = \pi'_j(u_1)$ and $v_j = \pi'_j(v_1)$. The rest of the argument is similar as in Case 1,
we just require that $\sigma_{1,x}(u_1) = u_x$ and $\sigma_{1,x}(v_1) = v_x$.

\emph{Case 3: The orbit $[e]$ is of size $\ell = {k \over 2}$.}
Then $e$ is projected to a half-edge in $H_{i+1}$.  The edge $e$ is halvable and corresponds to a
halvable proper atom or dipole $A$ in $G_i$.  Let $[e] = \{e_1,\dots,e_\ell\}$, and
$A_1,\dots,A_\ell$ be the corresponding atoms. Let $u_j$ be an arbitrary vertex of $e_j$ and let
$v_j$ be the second vertex of $e_j$, so $e_j = u_jv_j$.  Again, we arbitrarily choose $\bo$-isomorphisms
$\sigma_{1,x}$ from $A_1$ to $A_x$, for $x = 2,\dots,\ell$, such that $\sigma_{1,x}(u_1) = u_x$ and
$\sigma_{1,x}(v_1) = v_x$, and define $\sigma_{x,y} = \sigma_{1,y} \sigma^{-1}_{1,x}$.

Since $A_1$ is a halvable atom, there exists a semiregular involution $\tau_1 \in \Autbo{A_1}$ which exchanges
$u_1$ and $v_1$.  Then $\tau_1$ defines a semiregular involution of $A_x$ by conjugation as $\tau_x =
\sigma_{1,x} \tau_1 \sigma^{-1}_{1,x}$. It follows that
\begin{equation} \label{eq:tau}
\tau_y = \sigma_{x,y} \tau_x \sigma^{-1}_{x,y},\qquad\text{and consequently}\qquad \sigma_{x,y}
\tau_x = \tau_y \sigma_{x,y},\qquad \forall x,y \in \{1,\dots,\ell\}.
\end{equation}
We put $\hat\sigma_{x,y} = \sigma_{x,y} \tau_x = \tau_y \sigma_{x,y}$ which is a $\bo$-isomorphism
mapping $A_x$ to $A_y$ such that $\hat\sigma_{x,y}(u_x) = v_y$ and $\hat\sigma_{x,y}(v_x) = u_y$. We
note that $\hat\sigma_{x,x} = \tau_x$. In the extension, we put $\pi |_{\int A_x} = \sigma_{x,y}
|_{\int A_x}$ if $\pi'(u_x) = u_y$, and $\pi |_{\int A_x} = \hat\sigma_{x,y} |_{\int A_x}$ if
$\pi'(u_x) = v_y$. 

Aside (\ref{eq:sigma}), we get the following additional identities:
\begin{equation} \label{eq:hatsigma}
\hat\sigma_{x,z} = \sigma_{y,z} \hat\sigma_{x,y},\qquad
\hat\sigma_{x,z} = \hat \sigma_{y,z} \sigma_{x,y},\quad\text{and}\quad
\sigma_{x,z} = \hat \sigma_{y,z} \hat\sigma_{x,y},\qquad \forall x,y,z \in \{1,\dots,\ell\}.
\end{equation}
We just argue the last identity:
$$\hat \sigma_{y,z} \hat\sigma_{x,y} = \tau_z (\sigma_{y,z} \sigma_{x,y}) \tau_x = \tau_z
\sigma_{x,z} \tau_x = \tau_z \tau_z \sigma_{x,z} = \sigma_{x,z},$$
where the last equality holds since $\tau_z$ is an involution. It follows that for any $\pi'_1,
\pi'_2 \in \Gamma_{i+1}$ with extensions $\pi_1,\pi_2$, the composition $\pi_2 \pi_1$ is correctly
defined.

To conclude the proof, we argue semiregularity of $\Gamma_i$ as in Case 1. If $\pi \in \Gamma_i$
fixes a vertex $v$, then either $v$ is a vertex of $G_{i+1}$, or it is an internal vertex of one of
the atoms $A_x$.  In the second case, we get $\pi(A_x)=A_x$, so either $\pi |_{\int A_x} =
\sigma_{x,x} |_{\int A_x}=\id$, or $\pi |_{\int A_x} = (\tau_x\sigma_{x,x})|_{\int A_x}=\tau_x
|_{\int A_x}$.  Since $\tau_x$ is semiregular and $\pi$ fixes $v \in \int A_x$, $\pi$ cannot act as
$\tau_x$. In each case we conclude that $\pi$ is an extension of identity.  Then $\pi |_{\int A_y} =
\sigma_{y,y} |_{\int A_y} = \id$, for $y = 1,\dots,\ell$, so $\pi$ is the identity. If non-trivial
$\pi \in \Gamma_i$ fixes an edge $\hat e$, we need to argue that $\hat e$ is halvable. As in Case 1,
if $\hat e$ belongs to $G_{i+1}$, it is halvable. If $\hat e$ is an interior edge of $A_x$, then as
above $\pi(A_x) = A_x$. Therefore either $\pi$ is the identity, contradicting non-triviality of
$\pi$, or $\pi |_{\int A_x} = \tau_x |_{\int A_x}$ and since $\tau_x$ acts semiregularly, $\hat e$
is halvable.\qed
\end{proof}

\begin{corollary} \label{cor:all_group_expansions}
With the above notation, all possible extensions $\Gamma_i$ of $\Gamma_{i+1}$ are constructed using
the approach in the proof of Lemma~\ref{lem:group_expansion} by the following different choices for
each orbit $[e] = \{e_1,\dots,e_p\}$ of edges in $G_{i+1}$ corresponding to atoms $A_1,\dots,A_p$ of
$G_i$, respectively:
\begin{packed_itemize}
\item In Case 1, by all choices of $\sigma_{x,y} : A_x \to A_y$ such that $\sigma_{x,y}(u_x) = u_y$
and (\ref{eq:sigma}) holds.
\item In Case 2, by all choices of $\sigma_{x,y} : A_x \to A_y$ such that $\sigma_{x,y}(u_x) = u_y$,
$\sigma_{x,y}(v_x) = v_y$, and (\ref{eq:sigma}) holds.
\item In Case 3, by all choices of $\sigma_{x,y},\hat\sigma_{x,y} : A_x \to A_y$ and $\tau_x : A_x
\to A_x$ such that $\sigma_{x,y}(u_x) = u_y$, $\sigma_{x,y}(v_x) = v_y$, $\tau_x$ is a semiregular
involution, $\tau_x(u_x) = v_x$ and (\ref{eq:sigma}), (\ref{eq:tau}), and (\ref{eq:hatsigma}) hold.
\end{packed_itemize}
\end{corollary}

\begin{proof}
From the proof of Lemma~\ref{lem:group_expansion}, we know that for each choice of the input
parameters, we construct an extension $\Gamma_i$ of $\Gamma_{i+1}$. To show that all extensions are
covered, let $\Gamma_i$ be an arbitrary extension of $\Gamma_{i+1}$.  Let $|\Gamma_i| =
|\Gamma_{i+1}| = k = p$. Let $\pi \in \Gamma_i$ and $\pi' = \Phi_i(\pi)$.  If $\pi'(e_x) = e_y$,
then Lemma~\ref{lem:atom_automorphisms} and the definition of the reduction epimorphism $\Phi_i$
implies that $\pi|_{A_x}$ is a $\bo$-isomorphism from $A_x$ to $A_y$.

In Cases 1 and 2, we have $[e] = \{e_1,\dots,e_k\}$, corresponding to atoms $A_1,\dots,A_k$. From
semiregularity, for every $x,y$ there exists a unique isomorphism $\pi_{x,y} \in \Gamma_i$ such that
$\pi'_{x,y} = \Phi_i(\pi_{x,y})$ maps $e_x$ to $e_y$.  Let $\sigma_{x,y} = \pi_{x,y}|_{A_x}$. Since
$\Gamma_i$ forms a semiregular group, we get that
$$\pi_{y,z} \pi_{x,y} = \pi_{x,z},$$
so the property (\ref{eq:sigma}) is satisfied. The property $\sigma_{x,y}(u_x) = u_y$ (and in Case 2
also the property $\sigma_{x,y}(v_x) = v_y$) follows from the fact that vertices of $\bo
A_1,\dots,\bo A_k$ form one (in Case 1) or at most two (in Case 2) orbits of size $k$.

In Case 3, let $\ell = {k \over 2} = {p \over 2}$, $[e] = \{e_1,\dots,e_\ell\}$ and let $e_x$ consist of
darts $d_x$ and $\hat d_x$, incident with $u_x$ and $v_x$, respectively. Since $\Gamma_{i+1}$
is semiregular, the darts $\{d_1,\dots,d_\ell,\hat d_1, \dots \hat d_\ell\}$ and the vertices
$\{u_1,\dots,u_\ell,v_1,\dots,v_\ell\}$ form orbits of size $k$. From semiregularity, there exist
unique $\pi_{x,y}$ and $\hat\pi_{x,y}$ such that $\pi'_{x,y} = \Phi_i(\pi_{x,y})$ maps $d_x$ to
$d_y$ and $\hat\pi'_{x,y} = \Phi_i(\hat\pi_{x,y})$ maps $d_x$ to $\hat d_y$. Let $\sigma_{x,y} =
\pi_{x,y}|_{A_x}$, $\hat\sigma_{x,y} = \hat\pi_{x,y}|_{A_x}$ and $\tau_x = \hat\sigma_{x,x}$. Since
$\Gamma_i$ forms a semiregular group, we get that
$$\pi_{y,z} \pi_{x,y} = \hat\pi_{y,z} \hat\pi_{x,y} = \pi_{x,z}\qquad\text{and}\qquad
\pi_{y,z} \hat\pi_{x,y} = \hat\pi_{y,z} \pi_{x,y} = \hat\pi_{x,z}.$$
Therefore, the properties (\ref{eq:sigma}), (\ref{eq:tau}), and (\ref{eq:hatsigma}) hold.  Since
$\pi_{x,y}(u_x) = u_y$, $\pi_{x,y}(v_x) = v_y$, $\hat\pi_{x,y}(u_x) = v_y$, $\hat\pi_{x,y}(v_x) =
u_y$, the second required property follows.\qed
\end{proof}

As it was observed in the proof of Lemma~\ref{lem:group_expansion}, all these choices may be derived
from arbitrary choices of $\sigma_{1,2},\dots,\sigma_{1,k}$ in Cases 1 and 2 and from arbitrary choices
of $\sigma_{1,2},\dots,\sigma_{1,\ell},\tau_1$ in Case 3. The reader may notice that different
choices of $\sigma_{1,2},\dots,\sigma_{1,p}$ lead to isomorphic quotients $G_i / \Gamma_i$, but
different choices of $\tau_1$ might lead to non-isomorphic quotients.

In the language of commutative diagrams, given a group $\Gamma_{i+1}$,
Lemma~\ref{lem:group_expansion} proves that there exists an extension $\Gamma_i$ such that
Diagram~(\ref{eq:red_diagram}) commutes, while Corollary~\ref{cor:all_group_expansions} describes
all such extensions $\Gamma_i$.

There is a group theoretical reformulation of Lemma~\ref{lem:group_expansion} and of
Corollary~\ref{cor:all_group_expansions}. A pair $(\Gamma,D)$, where $\Gamma$ is a group acting on a
set $D$ is called a \emph{$\Gamma$-space}. A morphism between the spaces $(\Gamma,D)$ and
$(\Gamma',D')$ is a pair $(\Phi,f)$ such that $\Phi:\Gamma\to\Gamma'$ is a group homomorphism and
$f:D\to D'$ is a function, linked by the equation $f(g\cdot d)=\Phi(g)\cdot f(d)$. Let
$\textrm{red.}: \bD(G)\to \bD(G_r)$ be the mapping taking $d \mapsto d'$ if the full expansion of the dart $d'\in
\bD(G_r)$ contains $d$. The following corollary holds:

\begin{corollary} \label{thm:quotientprimitive}
Let $G = G_0,\dots,G_r$ be the reduction series, $\Phi=\Phi_{r-1} \circ
\Phi_{r-2} \circ \cdots \circ \Phi_0$, and $\Gamma\leq \Aut(G)$ be semiregular.
There is a uniquelly determined semiregular subgroup $\Gamma_r=\Phi(\Gamma)\cong\Gamma$ in
$\Aut(G_r)$ such that $( \Phi, \textrm{red.})$ is a morphism between the spaces $(\Gamma,\bD(G))$ and
$(\Gamma_r,\bD(G_r))$.  Moreover, Corollary~\ref{cor:all_group_expansions} describes all possible
ways of reconstruction of $(\Gamma,\bD(G))$ from $(\Gamma_r,\bD(G_r))$. 
\end{corollary}

\heading{Quotient Expansion.}
Recall the description of quotients of atoms from Section~\ref{sec:quotients_of_atoms}.  We are
ready to establish the main theorem of this paper proving the aforementioned meaning of quotient
expansion in Diagram~(\ref{eq:exp_diagram}). It states that every quotient $H_i$ of $G_i$ can be
constructed from some quotient $H_{i+1}$ of $G_{i+1}$ by replacing edges, loops and half-edges,
created from projections of edges $G_{i+1}$ replacing atoms of $G_i$, by corresponding edge-, loop-
and half-quotients.

\begin{proof}[Theorem~\ref{thm:quotient_expansion}]
Let $H_{i+1} = G_{i+1} / \Gamma_{i+1}$ and let $H_i$ be constructed in the above way. We first argue
that $H_i$ is a quotient of $G_i$, i.e., it is equal to $G_i / \Gamma_i$ for some $\Gamma_i$ extending
$\Gamma_{i+1}$. We use the construction from Corollary~\ref{cor:all_group_expansions}, where we
choose $\sigma_{1,2},\dots,\sigma_{1,p}$ arbitrarily and use the involutory semiregular automorphism
$\tau$ from the definition of half-projection as $\tau_1$. Let $k = |\Gamma_i| = |\Gamma_{i+1}|$.

It remains to argue that the constructed graph $G_i / \Gamma_i$ is isomorphic to $H_i$. Since
$\Gamma_i$ acts the same as $\Gamma_{i+1}$ outside interiors of atoms of $G_i$, we get that $G_i /
\Gamma_i$ and $H_i$ are isomorphic there. Next, we show that this partial isomorphism can be
extended to a full isomorphism. Let $[e] = \{e_1,\dots,e_p\}$ be an edge-orbit of $\Gamma_{i+1}$
and these edges correspond to atoms $A_1,\dots,A_p$. Let $p$ and $p'$ be the covering projections
$G_i \to G_i / \Gamma_i$ and $G_{i+1} \to G_{i+1} / \Gamma_{i+1}$, respectively.

If $p'(e)$ is a pendant edge, an edge, or a loop, we get Cases 1 and 2 in the proof of
Lemma~\ref{lem:group_expansion} and $p=k$. In $H_i$, we replace $p'(e)$ of $H_{i+1}$ by the edge- or
the loop-projection of, say, $A_1$. Similarly, in $G_i / \Gamma_i$, we have $p(\int A_1) \cong \int
A_1$ since the vertex- and edge-orbits of $\Gamma_i$ on $\int A_1$ are generated
by isomorphisms $\sigma_{1,2},\dots,\sigma_{1,n}$.

If $p'(e)$ is a half-edge, we get Case 3 in the proof of Lemma~\ref{lem:group_reduction} and $p =
\ell = {k \over 2}$. In $H_i$, we replace $p'(e)$ of $H_{i+1}$ by a half-quotient of, say, $A_1$,
constructed from a semiregular involutory automorphism $\tau_1$. The action of $\Gamma_i$ is
generated on atoms $A_1,\dots,A_\ell$ by $\sigma_{1,2},\dots,\sigma_{1,\ell},\tau_1$. Since $\tau_1$
is included, we get that $p(A_1)$ is isomorphic to the half-quotient. So $H_i$ and $G_i / \Gamma_i$
are isomorphic there as well.

On the other hand, if $H_i$ is a quotient, it replaces the edges, loops and half-edges of $H_{i+1}$
by some quotients, so we can generate $H_i$ in this way. The reason is that according to
Corollary~\ref{cor:all_group_expansions}, we can generate all $\Gamma_i$ extending the action of
$\Gamma_{i+1}$  onto the interiors of the atoms forming an orbit  by choosing
$\sigma_{1,2},\dots,\sigma_{1,p}$, and possibly $\tau_1$, and using the process repeatedly until the
action of $\Gamma_i$ is defined on the interiors of all atoms.\qed
\end{proof}

We say that two quotients $H_i$ and $H'_i$ extending $H_{i+1}$ are \emph{different} if there exists
no isomorphism of $H_i$ and $H'_i$ which fixes the vertices and edges common with $H_{i+1}$. (But
$H_i$ and $H'_i$ still might be isomorphic.) According to Lemma~\ref{lem:unique_quotients}, the edge
and loop-quotients are uniquely determined, so we are only free in choosing half-quotients.  For
non-isomorphic choices of half-quotients, we get different graphs $H_i$. For instance suppose that
$H_{i+1}$ contains a half-edge corresponding to the dipole from
Fig.~\ref{fig:exponentially_many_quotients}. Then in $H_i$ we can replace this half-edge by one of
the four possible half-quotients of this dipole.

\begin{corollary} \label{cor:unique_expansion}
If $H_{i+1}$ contains no half-edge, then $H_i$ is uniquely determined. Thus, for $\Gamma_r$ of an odd order, the quotient $H_r$ uniquely determines $H_0$.
\end{corollary}

\begin{proof}
This is implied by Theorem~\ref{thm:quotient_expansion} and Lemma~\ref{lem:unique_quotients} which
states that edge- and loop-quotients are uniquely determined. If the order of $\Gamma_r$ is odd, no
half-edges are constructed in $H_r$, so no half-quotients ever appear.\qed
\end{proof}

In the next two subsections we shall discuss some features of regular graph coverings in connection with the $3$-connected reduction that are important
from the algorithmic point of view. 

\subsection{Half-quotients of Dipoles}

In Lemma~\ref{lem:dipole_quotients_noncolored}, we describe that a dipole $A$ without colored edges can
have at most $\bigl\lfloor{\be(A) \over 2}\bigr\rfloor+1$ pairwise non-isomorphic half-quotients. This
statement can be easily altered to dipoles with colored edges which admit a much larger number of
half-quotients:

\begin{lemma} \label{lem:dipole_quotients}
Let $A$ be a dipole with colored edges. Then the number of pairwise non-isomorphic half-quotients is
bounded by $2^{\lfloor \be(A) / 2 \rfloor}$ and this bound is achieved. 
\end{lemma}

\begin{proof}
First, we analyze the structure of all involutory semiregular automorphisms $\tau$ acting on $\int
A$.  Concerning the non-halvable edges of $A$, the undirected edges of each color class have to be
paired by $\tau$ together. Further, each directed edge has to be paired with a directed edge of the
opposite direction and the same color. No matter how $\tau$ matches the non-halvable edges, if $m$
is a multiplicity in a colour class, in the half-quotient these edges are mapped onto $m/2$ loops. Hence
this part of the quotient is uniquely determined no matter which $\tau$ is used.

Now we discuss the action of $\tau$ on the remaining at most $\be(A)$ halvable edges of $A$. These
edges belong to $c$ color classes having $m_1,\dots,m_c$ edges.  Each automorphism $\tau$ has to
preserve the color classes, so it acts independently on each class.

We concentrate only on one color class having $m_i$ edges. Each edge $e$ of the class is either
fixed, or it is matched with another edge in the class.  In the first case $e$ projects to a
half-edge, in the second case the two edges swapped by $\tau$ are mapped onto a loop.  Denote by
$f(m_i)$ the number non-isomorphic half-quotients of the sub-dipole containing just halvable edges of
this class.  Then the total number of pairwise non-isomorphic half-quotients of $A$ is equal to
\begin{equation} \label{eq:dipole_bound}
\prod_{1 \le i \le c} f(m_i).
\end{equation}

\begin{figure}[t!]
\centering
\includegraphics{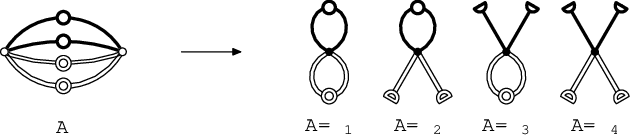}
\caption{An example of a dipole $A$ with a pair of black halvable edges and a pair of white halvable
edges. There exist four pairwise non-isomorphic half-quotients $A / \Gamma_1,\dots,A / \Gamma_4$
such that $\Gamma_1$ and $\Gamma_2$ swap the black edges, $\Gamma_3$ and $\Gamma_4$ fix the black
edges, $\Gamma_1$ and $\Gamma_3$ swap the white edges, and $\Gamma_2$ and $\Gamma_4$ fix the white
edges.}
\label{fig:exponentially_many_quotients}
\end{figure}

The rest of the proof is similar to the proof of Lemma~\ref{lem:dipole_quotients_noncolored}.  The
resulting half-quotient only depends on the number of fixed edges and fixed two-cycles in the
considered color class.  We can construct at most $f(m_i) = \lfloor{m_i \over 2}\rfloor + 1$
pairwise non-isomorphic half-quotients, since we may have zero to $\lfloor{m_i \over 2}\rfloor$ loops
with the complementing number of half-edges.

The value~(\ref{eq:dipole_bound}) is maximized when each class contains exactly two edges. (Except
for one class containing either three edges, or one edge if $\be(A)$ is odd.)

To prove that the upper bound is sharp, suppose first that the number of  the edges is even, and let
each color class contain two edges, say $i$-th class is formed by the edges $e_i$ and $e'_i$. For a
semiregular involution $\tau$ there are exactly two choices, either both $e_i$ and $e'_i$ are fixed
by $\tau$, or they form an orbit of size two in the action of $\tau$.  Since we can define $\tau$
for each color class independently, we get the required number of half-quotients employing all possible
involutions $\tau$. In the odd case, we assume that one color class contains exactly one edge, while the
orders are of size two. The single edge in the class must be fixed by $\tau$.  Further the
construction of both $\tau$ and of the half-quotient proceeds as in the even case.
Figure~\ref{fig:exponentially_many_quotients} shows an example for $\be(A) = 4$.
\qed
\end{proof}

Assume that $H_{i+1}$ contains a half-edge corresponding to a half-quotient of a dipole in $H_i$.
By Theorem~\ref{thm:quotient_expansion}, the number of non-isomorphic expansions $H_i$ of $H_{i+1}$
can be exponential in the size difference of $H_i$ and $H_{i+1}$.

\subsection{The Block Structure of Quotients}

We show how the block structure is preserved during expansions. A block atom $A$ of $G_i$ is always
projected by an edge-projection, and so it corresponds to a block atom of $H_i$. It remains to deal
with a proper atom or a dipole $A$, and let $\bo A = \{u,v\}$.

For an edge-projection $p |_A$, we get $p(u) \ne p(v)$, and $p(A)$ is isomorphic to an atom in
$H_i$.

For a loop- or a half-projection $p |_A$, we get $p(u) = p(v)$ and $p(u)$ is an articulation
of $H_i$. If $A$ is a dipole, then $p(A)$ is a pendant star of half-edges and loops attached to
$p(u)$. If $A$ is a proper atom, we use the characterization by Lemma~\ref{lem:proper_atoms}.
\begin{packed_itemize}
\item \emph{For a loop-projection $p |_A$,} we have $p(A)$ either essentially a cycle (when $A^+$ is
essentially a cycle), or a pendant block with attached single pendant edges (when $A^+$ is
essentially 3-connected).
\item \emph{For a half-projection $p |_A$,} we have $p(A)$ either a path ending with a half-edge and
with attached single pendant edges (when $A^+$ is essentially a cycle), or a pendant block with
attached single pendant edges and half-edges (when $A^+$ is essentially 3-connected).  In the last
case, the only articulations other than $p(u) = p(v)$ separate single pendant edges and half-edges,
since the fiber over an articulation in a 2-fold cover is a 2-cut, so $A$ would not be a proper atom
otherwise.
\end{packed_itemize}

\begin{lemma} \label{lem:block_structure}
The block structure of $H_{i+1}$ is preserved in $H_i$, possibly with some new subtrees of blocks
attached.
\end{lemma}

\begin{proof}
By Theorem~\ref{thm:quotient_expansion}, edges inside blocks are replaced by edge-quotients of block
atoms, proper atoms and dipoles which preserves 2-connectivity. New rooted subtrees of blocks in $H_i$ are
created by replacing pendant edges with block atoms, loops by loop-quotients, and half-edges by
half-quotients.\qed
\end{proof}
 
\section{Planar Graphs} \label{sec:planar_graphs}

In this section, we show implications of our theory to planar graphs. We first discuss some well-known
properties of automorphism groups of 3-connected planar graphs. We use them to characterize the
quotients of planar graphs which results in a proof of Negami's Theorem.  The key point is
that semiregular groups of automorphisms of 3-connected planar graphs
are well understood, and therefore one can  effectively determine
all regular covering projections  defined on a 3-connected planar graph.

\subsection{Automorphism Groups of 3-connected Planar Graphs}

We review  geometric properties of automorphism groups of 3-connected planar graphs.  By Whitney
Theorem~\cite{whitney1932congruent}, 3-connected planar graphs have unique embeddings onto the
sphere up to a homeomorphism. This is strengthen by Mani Theorem~\cite{mani1971automorphismen}
stating that the unique embedding can be realised on the sphere such that  all automorphisms of the
graph correspond to isometries of the underlying sphere.  Using these properties, we describe
possible automorphism groups of planar atoms and primitive graphs.

\begin{figure}[t!]
\centering
\includegraphics{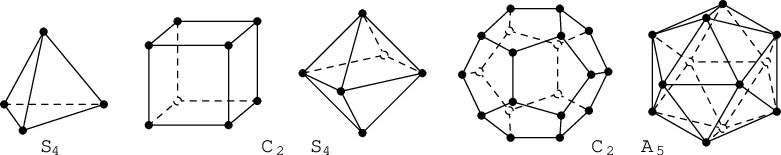}
\caption{The five platonic solids together with their automorphism groups.}
\label{fig:platonic_solids}
\end{figure}

\heading{Spherical Groups.} A group is \emph{spherical} if it is the group of the isometries of a
(finite) tiling of the sphere. The first class of spherical groups are the subgroups of the automorphism
groups of the platonic solids. Their automorphism groups are isomorphic to $\gS_4$ for the tetrahedron, $\gS_4
\times \gC_2$ for the cube and the octahedron, and $\gA_5 \times \gC_2$ for the dodecahedron and the
icosahedron; see Fig.~\ref{fig:platonic_solids}.  Table~\ref{tab:number_of_subgroups} shows the
number of conjugacy classes of subgroups of these three groups.  Note that conjugate semiregular
subgroups $\Gamma$ determine isomorphic quotients $G / \Gamma$ of the one-skeletons.  The second class of spherical
groups is formed by four infinite families $\gC_n$, $\gD_n$, $\gC_n \times \gC_2$, and $\gD_n \times
\gC_2$, $n\geq 2$. All they act as groups of automorphisms of $n$-sided prisms.

\begin{table}[b!]
\centering
\small{
\begin{tabular}{|c|c|c|c|}
\hline
\multicolumn{4}{|c|}{$\gS_4$ of the order 24}\\
\hline
Order & Number & Order & Number\\
\hline
1 & 1 & 6 & 1\\
2 & 2 & 8 & 1\\
3 & 1 & 12 & 1\\
4 & 3 & & \\
\hline
\end{tabular}
\hskip 0.3cm
\begin{tabular}{|c|c|c|c|}
\hline
\multicolumn{4}{|c|}{$\gC_2 \times \gS_4$ of the order 48}\\
\hline
Order & Number & Order & Number\\
\hline
1 & 1 & 8 & 7\\
2 & 5 & 12 & 2\\
3 & 1 & 16 & 1\\
4 & 9 & 24 & 3\\
6 & 3 & & \\
\hline
\end{tabular}\\[0.3cm]
\begin{tabular}{|c|c|c|c|}
\hline
\multicolumn{4}{|c|}{$\gC_2 \times \gA_5$ of the order 120}\\
\hline
Order & Number & Order & Number\\
\hline
1 & 1 & 8 & 1\\
2 & 3 & 10 & 3\\
3 & 1 & 12 & 2\\
4 & 3 & 20 & 1\\
5 & 1 & 24 & 1\\
6 & 3 & 60 & 1\\
\hline
\end{tabular}}
\caption{The number of conjugacy classes of the subgroups of the groups of platonic solids.}
\label{tab:number_of_subgroups}
\end{table}

%\heading{Maps.}
%A (spherical) map $\calM$ is a 2-cell decomposition of the sphere $S$. A map is usually defined by a
%2-cell embedding of a connected graph $i \colon G \hookrightarrow S$. The connected components of $S
%\setminus i(G)$ are called faces of $\calM$.  An automorphism of a map is an automorphism of the
%graph preserving the incidences between the vertices, edges, and faces. Clearly, $\Aut(\calM)$ is one of
%the spherical groups and with the exception of paths and cycles, it is a subgroup of $\Aut(G)$. As a
%consequence of Whitney's theorem~\cite{whitney} we have the following.
%
%\begin{theorem}\label{lem:aut_3_conn}
%Let $\calM$ be the map given by the unique $2$-cell embedding of a $3$-connected graph into the
%sphere. Then $\Aut(G) \cong \Aut(M)$.\qed
%\end{theorem}

\heading{Isometries of the Sphere.}
We recall some basic definitions from geometry of the sphere; see~\cite[Chapter 3]{stillwell} and \cite[Section
6.II]{vca}.  An isometry is a distance-preserving homeomorphism of the unit sphere. Since isometries
are closed under composition, they form a group, denoted $\Iso(S^2)$. A first example of an isometry
is a \emph{rotation} with the axis passing through two opposite points of the sphere. A second
example is a \emph{reflection} by the plane passing through two opposite points of the sphere. A
third example is the \emph{antipodal mapping} which maps each point of the sphere to the opposite
point.  Each rotation fixes exactly two points that are opposite on the sphere, while each reflection fixes a great circle
passing through a pair of opposite points,  the antipodal mapping fixes no point of the sphere.  For the
purpose of this paper, it is sufficient to understand these isometries geometrically. They can be
defined algebraically, for instance as certain $3 \times 3$ real orthogonal matrices~\cite[Chapter
3]{strang_la}, or as certain M\"obius transformations~\cite[Chapter 3]{vca}.

In~\cite[Section 6.II]{vca}, it is proved that all isometries of the sphere are generated by all
rotations and reflections of the sphere. An isometry is called \emph{orientation preserving} (direct
or conformal in~\cite[Section 6.II]{vca}) when it preserves a chosen orientation of the sphere. It is
called \emph{orientation reversing} (called opposite or anticonformal in~\cite[Section 6.II]{vca})
if it changes the orientation of the sphere. It is proved in~\cite[Section 6.II]{vca} that every
orientation preserving isometry of the sphere is a rotation, while every orientation reversing
isometry is the composition of a rotation and a reflection.  For instance the antipodal mapping is
the orientation reversing isometry formed as the composition of a $180^\circ$ rotation with the
reflection whose defining plane is perpendicular to the axis of this rotation.

A subgroup of a group of isometries is called \emph{orientation preserving} if it contains only
orientation preserving isometries, and \emph{orientation reversing} otherwise.  Every orientation
reversing subgroup contains an orientation preserving subgroup of index two since the composition of
two orientation reversing automorphisms is an orientation preserving automorphism.

Precisely the following isometries of the sphere are involutions: the identity, every reflection,
every $180^\circ$ rotation, and the antipodal mapping. 

\heading{Geometry of Automorphisms.}
Isometries of $\R^3$ are distance preserving mappings of homeomorphisms of $\R^3$, and they form the
group $\Iso(\R^3)$. For a polyhedron $P$ embedded in $\R^3$, we denote $\bV(P)$ its vertices,
$\bE(P)$ its edges, and $\Iso(P)$ the subgroup of $\Iso(\R^3)$ of all isometries preserving $P$.
Mani~\cite{mani1971automorphismen} gives the following insight into geometry of automorphism groups
of 3-connected graphs.

\begin{theorem}[Mani~\cite{mani1971automorphismen}] \label{thm:mani}
Let $G$ be a 3-connected planar graph. There exists an associated polyhedron $P$ with $\bV(P) =
\bV(G)$ and $\bE(P) = \bE(G)$ such that $\Aut(G) \cong \Iso(P)$, i.e., every automorphism $\pi \in
\Aut(G)$ induces some isometry $\hat\pi \in \Iso(P)$ and vice versa.
\end{theorem}

In other words, there exists a polyhedron $P$ in $\R^3$ with vertices corresponding to $\bV(G)$ and
edges corresponding to $\bE(G)$ such that each automorphism $\pi \in \Aut(G)$ induces an isometry of
$P$ and vice versa. Figure~\ref{fig:platonic_solids} gives examples of such polyhedra associated to
the graphs of platonic solids. 

Notice that every isometry in $\Iso(P)$ preserves the center of $P$.  Suppose that we place the
polyhedron $P$ into $S^2$ so that the centers of $P$ and $S^2$ coincide, and project $P$ onto $S^2$.
Then each isometry of $P$ corresponds to some isometry of the sphere and $\Iso(\R^3) > \Iso(S^2) >
\Iso(P)$. If we view the projection of $P$ onto $S^2$ as an embedding of $G$, called a
\emph{geometric embedding}, every automorphism of $G$ corresponds to some isometry of the sphere,
i.e., to a rotation or to a composition of a rotation and a reflection. Therefore, the
aforementioned geometric name of isometries translate to automorphisms of 3-connected planar graphs,
so for instance, an automorphism of a 3-connected planar graph $G$ is called a rotation if the
corresponding isometry of the sphere is a rotation. Namely, it follows that $\Aut(G)$ is isomorphic
to one of spherical groups.

\heading{Stabilizers.}
Let $G$ be a 3-connected planar graph and $u \in \bV(G)$. The stabilizer of $u$ in $\Aut(G)$ is a
subgroup of a dihedral group and it has the following description in the language of isometries. If
$\Stab(u) \cong \gC_n$, for $n \ge 3$, it is generated by a rotation of order $n$ that fixes $u$ and
the opposite point of the sphere. The opposite point of the sphere may be another vertex or a center
of a face.  If $\Stab(u) \cong \gD_n$, for $n \ge 2$, it consists of rotations fixing $u$ and the
opposite point of the sphere and reflections fixing a great circle passing through $u$ and the
opposite point. Each reflection always fixes a great circle, containing aside $u$ at least two other
points of the geometric embedding of $G$, each being either a center of some edge, or another
vertex. When $\Stab(u) \cong \gD_1 \cong \gC_2$, it is generated either by a $180^\circ$ rotation
or by a reflection; in the former case, the opposite point of the sphere may also be the center of
an edge.

Let $e \in \bE(G)$ and $e = uv$. The stabilizer of $e$ in $\Aut(G)$ is a subgroup of $\gC_2^2$. When $\Stab(e)
\cong \gC_2^2$, it contains the following three non-trivial isometries. First, the $180^\circ$
rotation fixing the center of $e$ and the opposite point of the sphere that is a vertex, the center of
an edge, or the center of an even face. Next, two reflections perpendicular to each other fixing
great circles passing through the center of $e$, one fixing both $u$ and $v$, the other swapping
them. When $\Stab(e) \cong \gC_2$, it is generated by only one of these three isometries.

\subsection{Automorphism Groups of Planar Primitive Graphs and Atoms}

Mani Theorem~\ref{thm:mani} allows to describe possible automorphism groups of planar atoms and
primitive graphs which appear in the reduction tree for a planar graph $G$. First, we describe the
automorphism groups of planar primitive graphs.

\begin{lemma} \label{lem:planar_primitive_graph}
The automorphism group $\Aut(G)$ of a planar primitive graph $G$ is a spherical group.
\end{lemma}

\begin{proof}
Recall that a graph is essentially 3-connected if it is a 3-connected graph with attached single
pendant edges to some of its vertices. If $G$ is essentially 3-connected, then $\Aut(G)$ is a
spherical group by Mani Theorem~\ref{thm:mani}, Lemma~\ref{lem:aut_ess_3-conn}, and the fact that
spherical groups is closed under taking subgroups. If $G$ is $K_1$, $K_2$ or $C_n$ with
attached single pendant edges, then $\Aut(G)$ is, respectively, trivial, a subgroup of $\gC_2$, or
of $\gD_n$.\qed
\end{proof}

\begin{figure}[t!]
\centering
\includegraphics{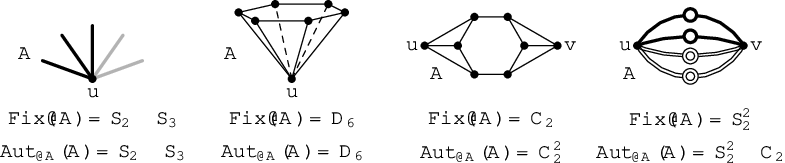}
\caption{An atom $A$ together with its groups $\Fix(\bo A)$ and $\Autbo{A}$. From left to right, a
star block atom, a non-star block atom, a proper atom, and a dipole.}
\label{fig:automorphisms_of_atoms}
\end{figure}

Next, we deal with the automorphism groups of planar atoms.  Let $A$ be a planar atom. Recall that
$\Autbo{A}$ is the set-wise stabilizer of $\bo A$, and $\Fix(\bo A)$ is the point-wise stabilizer of
$\bo A$. The following lemma determines $\Autbo{A}$ and $\Fix(\bo A)$; see
Fig.~\ref{fig:automorphisms_of_atoms} for examples.

\begin{lemma} \label{lem:planar_atom_aut_groups}
Let $A$ be a planar atom.
\begin{packed_enum}
\item[(a)] If $A$ is a star block atom, then $\Autbo{A} = \Fix(\bo A)$ which is a direct product of
symmetric groups. 
\item[(b)] If $A$ is a non-star block atom, then $\Autbo{A} = \Fix(\bo A)$ and it is a subgroup of
a dihedral group.
\item[(c)] If $A$ is a proper atom, then $\Autbo{A}$ is a subgroup of $\gC_2^2$ and $\Fix(\bo A)$ is
a subgroup of $\gC_2$.
\item[(d)] If $A$ is a dipole, then $\Fix(\bo A)$ is a direct product of symmetric groups. If $A$ is
symmetric or halvable, then $\Autbo{A} = \Fix(\bo A) \times \gC_2$. If $A$ is asymmetric, then
$\Autbo{A} = \Fix(\bo A)$.
\end{packed_enum}
\end{lemma}

\begin{proof}
(a) Since $|\bo A|=1$, we have $\Autbo{A} = \Fix(\bo A)$.  Since the edges of each color class of
the star block atom $A$ can be arbitrarily (and independently) permuted, $\Fix(\bo A)$ is a direct
product of symmetric groups.

(b) Similarly as in Case (a), $|\bo A|=1$, and we have $\Autbo{A} = \Fix(\bo A)$.  We construct a graph $B$ from $A$
by removing all single pendant edges. By Lemma~\ref{lem:aut_ess_3-conn}, $\Autbo{A} \le \Autbo{B}$. By
Lemma~\ref{lem:non_star_block_atoms}, either $B$ is a cycle, $K_2$, or a 3-connected
planar graph. In the first two cases, $\Autbo{B}$ is a subgroup of $\gC_2$, while in the last case, it is the
stabilizer of a vertex in a 3-connected planar graph which is a subgroup of $\gD_n$, where $n$ is the degree
of the articulation separating $A$ in the subgraph $B$.

(c) Let $A$ be a proper atom with $\bo A = \{u,v\}$. As in the previous case, we construct the
subgraph $B$ from $A$ by removing all the single pendant edges, and let $B^+ = B + uv$. Then
$\Autbo{A^+} \le \Autbo{B^+}$ and $\Fix(\bo A^+) \le \Fix(\bo B^+)$ by
Lemma~\ref{lem:aut_ess_3-conn}. We have that $\Autbo{B}$ is isomorphic to the stabilizer of $uv$ in
$\Aut(B^+)$ and $\Fix(\bo B^+) \cong \Fix(\bo B)$, where $B^+=B+ uv$.  By
Lemma~\ref{lem:proper_atoms}, $B^+$ is either a cycle, or a 3-connected planar graph. In the former
case, $\Aut(B^+)$ is a subgroup of $\gC_2$ and $\Fix(\bo B^+)$ is trivial. In the latter case, the
stabilizer of $uv$ in $\Aut(B^+)$ is a subgroup of $\gC_2^2$ while $\Fix(\bo B^+)$ is a subgroup of
$\gC_2$.

(d) For an asymmetric dipole, we have $\Autbo{A} = \Fix(\bo A)$ which is a direct product of
symmetric groups. For a symmetric or halvable dipole, we can swap the vertices in $\bo A$ by an
involution $\tau$ fixing all the edges.  Since $\tau$ commutes with all the elements in $\Fix(\bo
A)$ we have that $\Autbo{A} \leq \Fix(\bo A)\times\langle\tau\rangle$. On the other hand, the
product of any two automorphisms swapping the two vertices belongs to $\Fix(\bo A)$, hence $\Fix(\bo
A)$ is a normal subgroup of index two. Thus $\Autbo{A} \cong \Fix(\bo A)\times\langle\tau\rangle$. \qed
\end{proof}

\subsection{Quotients of Planar Graphs and Negami's Theorem} \label{sec:quotients_of_planar_graphs}

In this section, we describe quotients of planar graphs geometrically. Using
Theorem~\ref{thm:quotient_expansion}, it only remains to understand the quotients of planar
primitive graphs and the half-quotients of planar proper atoms. We also show that our structural
theory gives a proof of Negami's Theorem~\cite{negami} as a straightforward by-product.

\heading{Quotients of the Sphere.}
The (branched) regular covering projections of surfaces are defined analogously as of graphs, where
the triviality of the stabilizers  of the defining group of isometries is required with the
exception of singular points called \emph{branch points}; see~\cite[Section 2.VI]{vca}. The
following well-known statement,  characterizing the regular quotients of the sphere, is a
consequence of the Riemann-Hurwitz equation:

\begin{lemma} \label{lem:quotients_of_sphere}
Let $\Gamma$ be a semiregular subgroup of isometries of the sphere $S^2$.
\begin{packed_enum}
\item[(a)] When $\Gamma$ is orientation preserving, then $S^2 / \Gamma$ is homeomorphic to the
sphere.
\item[(b)] When $\Gamma$ is orientation reversing and does not contain the antipodal mapping, then
$S^2 / \Gamma$ is homeomorphic to the disk.
\item[(c)] When $\Gamma$ is orientation reversing and contains the antipodal mapping, then $S^2 /
\Gamma$ is homeomorphic to the projective plane.
\end{packed_enum}
\end{lemma}

Correctness of this lemma may be argued directly as follows. If $\Gamma$ has two semiregular complementary
subgroups $\Gamma_1, \Gamma_2$ (i.e., $\Gamma_1,\Gamma_2 \le \Gamma$ and $\langle \Gamma_1 \cup
\Gamma_2 \rangle = \Gamma$), then
$$S^2 / \Gamma \cong (S^2 / \Gamma_1) / \Gamma_2.$$
Using this, it is sufficient for (a) to understand what are the quotients $S^2 / \langle \rho
\rangle$ where $\rho$ is a rotation of order $n$, so $\langle \rho \rangle \cong \gC_n$. Each
point-orbit consists of one point and its images obtained by repeatedly rotating this point by the
angle $360^\circ/n$. The quotients $S^2 / \langle \rho \rangle$ is one spherical luna of dihedral
angle $360^\circ/n$ with boundaries glue together which is homeomorphic to $S^2$.

To establish (b) and (c), each orientation reversing group $\Gamma$ contains an orientation
preserving subgroup $\Gamma'$ of index two. It is possible to write $\Gamma = \bigl\langle \Gamma'
\cup \{\tau\}\bigr\rangle$, where $\tau$ is a semiregular involution which is a reflection in (b),
and the antipodal mapping in (c). From (a), $S^2 / \Gamma'$ is homeomorphic to sphere. It is easy to
observe $S^2 / \langle \tau \rangle$ is the disk when $\tau$ is a reflection, and is the projective
plane when $\tau$ is the antipodal mapping.  

\begin{figure}[b!]
\centering
\includegraphics{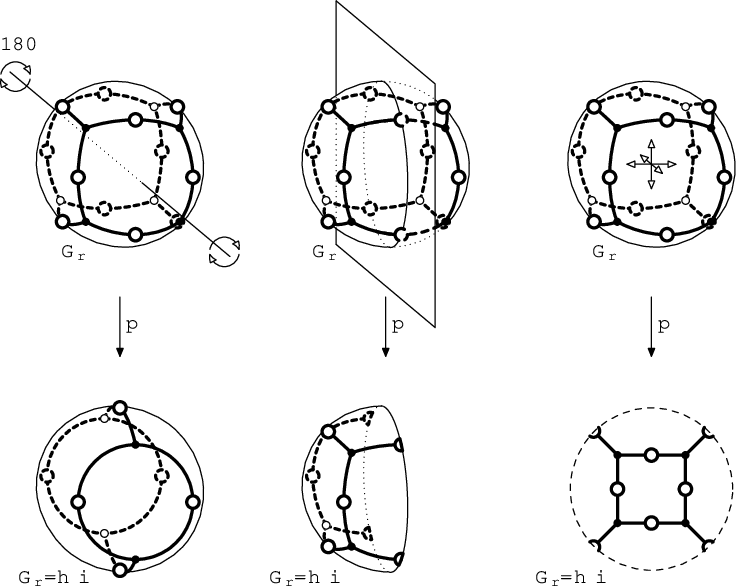}
\caption{From left to right, a rotational quotient, a reflectional quotient and an antipodal
quotient of the cube; also see Fig.~\ref{fig:quotients_of_cube}.}
\label{fig:quotients_of_primitive_graphs}
\end{figure}

For more details, we refer the reader to characterization of spherical groups and orbifolds (which
are regular quotients of the sphere, together with additional information about branch points)
in~\cite{conway2016symmetries}. See Fig.~\ref{fig:quotients_of_primitive_graphs} for examples of
these quotients of the sphere.

\heading{Quotients of 3-connected Planar Graphs.}
Using Lemma~\ref{lem:quotients_of_sphere}, we may characterize all possible quotients of 3-connected
planar graphs:

\begin{lemma} \label{lem:quotients_of_planar_graphs}
Let $G$ be a 3-connected planar graph and $\Gamma$ be a semiregular subgroup of $\Aut(G)$.
There are three types of quotients of $G$:
\begin{packed_head_enum}{(a)}
\item[(a)] \emph{Rotational quotients} -- The action of $\Gamma$ is orientation preserving and the
quotient $G / \Gamma$ is planar.
\item[(b)] \emph{Reflectional quotients} -- The action of $\Gamma$ is orientation reversing and does
not contain the antipodal mapping. Then the quotient $G / \Gamma$ is planar and necessarily
contains at least one half-edge. Further for $|\Gamma| = 2$, the quotient $G / \Gamma$ contains at
least three half-edges.
\item[(c)] \emph{Antipodal quotients} -- The action of $\Gamma$ is orientation reversing and
contains the antipodal mapping. Then $G / \Gamma$ is projective planar.
\end{packed_head_enum}
\end{lemma}

\begin{proof}
By Mani Theorem~\ref{thm:mani}, there exists a geometric embedding of $G$ into the sphere such that
each automorphism of $G$ extends to an isometry of the sphere. Let $\widetilde\Gamma$ be the
semiregular subgroup of isometries of the sphere corresponding to $\Gamma$. By
Lemma~\ref{lem:quotients_of_sphere}, we get that $S^2 / \widetilde\Gamma$ is the sphere (a), the
disk  (b), or the projective plane (c). The key observation is that if $G\hookrightarrow S^2$ is an
embedding, then it induces the embedding $G /\Gamma\hookrightarrow S^2 / \widetilde\Gamma$ of the
quotient graph into the quotient surface. So $G / \Gamma$ is planar for (a) and (b), and projective
planar for (c). 

Half-edges are created in $G / \Gamma$ when the embedding of $G$ places centers of edges into fixed
points of non-trivial isometries in $\widetilde\Gamma$. Each rotation only fixes two opposite points of the
sphere, so it might create half-edges only when its order is two. But a reflection $\tau$ in (b) fixes a
great circle of $S^2$. Since the action of $\tau$ is semiregular, no vertex can be placed on this
circle. Therefore, when the embedding of $G$ crosses this great circle from one hemisphere to the
other, the center of an edge is placed on this great circle. Since the embedding of $G$ is symmetric
with respect to $\tau$ along this great circle and $G$ is 3-connected, there are at least
three edges crossing the great circle, resulting into at least three half-edges in $G / \langle \tau
\rangle$. At least one of these half-edges remains in $G / \Gamma$.\qed
\end{proof}

Figure~\ref{fig:quotients_of_primitive_graphs} shows examples of these types of quotients. We note
that the antipodal quotient of a planar graph may, or may not, be planar; for an example, see
Fig.~\ref{fig:big_picture}.

\heading{Quotients of Primitive Graphs.} By Lemma~\ref{lem:primitive_graphs}, we know that every
primitive graph $G_r$ is either 3-connected with attached single pendant edges, or $K_2$ or $C_n$
with attached single pendant edges. These attached single pendant edges only make $\Aut(G_r)$
smaller by Lemma~\ref{lem:aut_ess_3-conn}, and correspond to single pendant edges in $G_r / \Gamma$.
Therefore it is sufficient to understand how possible quotients can look for 3-connected planar
graphs, $K_2$ and $C_n$. 

The quotients of 3-connected planar graphs are described in	Lemma~\ref{lem:quotients_of_planar_graphs}.
The quotients of $K_2$ are straightforward. Next, we characterize quotients of cycles, which completes
the description of possible quotients of primitive graphs:

\begin{lemma} \label{lem:cycle_quotients}
Let $\Gamma$ be a semiregular subgroup of $\Aut(C_n)$. Then $C_n / \Gamma$ is either a cycle, or a
path with two half-edges attached to its ends (only for $n$ even).
\end{lemma}

\begin{proof}
The former case happens when $\Gamma$ if generated by a rotation, i.e., $\Gamma \cong \gC_k$ for
some $k$. The latter case happens when $\Gamma$ contains a reflection, necessarily fixing the
centers of two edges, corresponding to two half-edges in $C_n / \Gamma$. In the latter case,
$\Aut(C_n) \cong \gD_k$ for some even $k$, otherwise each reflection fixes one edge and one vertex,
so its action is not semiregular.\qed
\end{proof}

\begin{figure}[b!]
\centering
\includegraphics{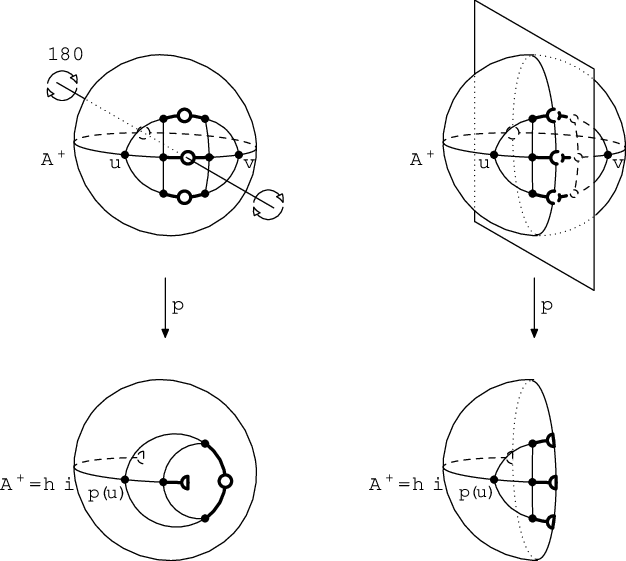}
\caption{The rotational quotient and reflectional quotient of a planar proper atom $A$ with the added edge
$uv$.}
\label{fig:halfquotients_of_proper_atoms}
\end{figure}

\heading{Half-quotients of Proper Atoms.} Next, we characterize the half-quotients of planar proper
atoms. For a proper atom $A$ with $\bo A = \{u,v\}$, we characterize all semiregular involutions
$\tau \in \Autbo{A}$ exchanging $u$ and $v$. For the extended proper atom $A^+=A+uv$, we get that
$\tau$ corresponds to a semiregular involution in $\Aut(A^+)$ fixing the added edge $uv$.  By
Lemma~\ref{lem:proper_atoms}, $A^+$ is either essentially a cycle, or essentially 3-connected.

When $A^+$ is essentially a cycle, $\tau$ can only be a reflection through $uv$ and another edge.
Therefore, the half-quoutient $A / \langle \tau \rangle$ is a path ending with a half-edge attached
to the last vertex and with single pendant edges attached some of the vertices with the exception of
the first one.

When $A^+$ is essentially 3-connected, we get the following two types of half-quotients (see
Fig.~\ref{fig:halfquotients_of_proper_atoms} for examples):

\begin{lemma} \label{lem:planar_proper_half_quotients}
Let $A$ be a planar proper atom such that $A^+$ is essentially 3-connected and let $\bo A =
\{u,v\}$. There are at most two half-quotients $A / \left<\tau\right>$ where $\tau \in \Autbo{A}$ is
an involutory semiregular automorphism transposing $u$ and $v$:
\begin{packed_enum}
\item[(a)] \emph{The rotational half-quotient} -- The involution $\tau$ is orientation preserving
and $A / \left<\tau\right>$ is planar with at most one half-edge.
\item[(b)] \emph{The reflectional half-quotient} -- The involution $\tau$ is a reflection and $A /
\left<\tau\right>$ is planar with at least two half-edges.
\end{packed_enum}
\end{lemma}

\begin{proof}
The graph $A^+$ is an essentially 3-connected planar graph, and let $B^+$ be constructed from $A^+$
by removal of all the pendant edges. It is sufficient to understand half-quotients of $B$, since
adding single pendant edges only restricts possible symmetries of $A^+$. Using Mani
Theorem~\ref{thm:mani}, consider a geometric embedding of $B^+$ into the sphere.  The stabilizer of
the edge $uv$ in $B^+$ is a subgroup of $\gC_2^2$, so there are at most two involutions $\tau$
fixing $uv$ and exchanging $u$ with $v$.

First, $\tau$ might be the $180^\circ$ rotation with the axis passing through the center of $uv$ and the
opposite point of the sphere. Then $A / \langle \tau \rangle$ contains at most one half-edge, when
an edge of $B^+$ is placed to this opposite point of the sphere fixed by $\tau$. 

Second, $\tau$ might be a reflection by the plane passing through the center of $uv$, perpendicular
to the segment $uv$. Simiraly, as in the proof of Lemma~\ref{lem:quotients_of_planar_graphs}, we argue that by $3$-connectivity there are
at least two edges other than $uv$ whose centers belong to this plane, resulting into at least two
half-edges in $A / \langle \tau \rangle$.

The semiregular subgroup $\langle \tau \rangle$ is either of type (a), or of type (b) in
Lemma~\ref{lem:quotients_of_planar_graphs}. Therefore the half-quotients $A /
\langle \tau \rangle$ is planar.\qed
\end{proof}

Suppose a planar proper atom $A$ has half-quotients of both  types (a) and (b). Since they have different numbers of half-edges,  they are
necessarily non-isomorphic.

\heading{Proof of Negami's Theorem.}
Using the above statements, we give a proof of Negami's Theorem~\cite{negami}. This theorem states
that a graph $H$ has a finite planar regular cover $G$ (i.e, $G/\Gamma \cong H$ for some semiregular
$\Gamma \le \Aut(G)$), if and only if $H$ is projective planar. For a given projective planar graph
$H$, the construction of a planar graph $G$ covering $H$ is easy: by embedding $H$ into the projective plane and
taking the antipodal double cover of this embedding, we get the graph $G$ embedded to the sphere.  Below, we
prove the harder implication:

\begin{theorem}[Negami~\cite{negami}] \label{thm:negami}
Let $G$ be a planar graph. Then every (regular) quotient of $G$ is projective planar.
\end{theorem}

\begin{proof}
We apply the reduction series on $G$ which produces graphs $G = G_0,G_1,\dots,G_r$ such that $G_r$
is primitive. If $G_r$ is essentially 3-connected, then by
Lemma~\ref{lem:quotients_of_planar_graphs} every quotient $H_r = G_r / \Gamma_r$ is projective
planar. If $G_r$ is $K_2$ or $C_n$ with single pendant edges attached, then by
Lemma~\ref{lem:cycle_quotients} every quotient $H_r = G_r / \Gamma_r$ is even planar.

By Theorem~\ref{thm:quotient_expansion}, every quotient $H = G / \Gamma$ can be constructed from
some $H_r$ by an expansion series in which we replace edges, loops and half-edges by edge-quotients,
loop-quotients and half-quotients, respectively. All edge- and loop-quotients are clearly planar. By
Lemma~\ref{lem:planar_proper_half_quotients}, every half-quotient of a proper proper atom is planar, and by
Lemma~\ref{lem:dipole_quotients} every half-quotient of a dipole is a set of loops and half-edges
attached to a single vertex, which is also planar. Therefore, these replacements can be done in a
way that the underlying surface of $H_r$ is not changed, so $H$ is also projective planar.\qed
\end{proof}

First note that since our definition of a regular covering is more general than the one used by
Negami,  we have proved a more general statement.  Further, our method gives a deeper insight into
the structure of the regular quotients of planar graph. A brief discussion follows.

Deciding whether $H = G / \Gamma$ is planar or non-planar projective is done on the primitive graph
$G_r$. It is non-planar if and only if $\Gamma_r$ contains a semiregular antipodal involution and
the resulting quotient $H_r = G_r / \Gamma_r$ is non-planar. A more precise analysis when $H_r$ is
planar and when not would be possible by checking the action of all possible spherical groups and
the resulting quotients.  For instance, the cube in Fig.~\ref{fig:quotients_of_cube} has as the
half-quotient $K_4$ induced by the antipodal mapping which is in fact planar. On the other hand,
the half-quotient of the dodecahedron generated by the antipodal mapping is the famous Petersen
graph which is non-planar.

\heading{Comparison with Negami's Proof.} We describe Negami's approach~\cite{negami}, adapted
to our notation. First, if $G$ is a 3-connected planar graph, a similar argument~\cite[Section
2]{negami} as in the proof of Lemma~\ref{lem:quotients_of_planar_graphs} is used to show that $G /
\Gamma$ is either planar, or projectively planar. So it remains to deal with the situation that $G$
is not 3-connected.

Let $G$ be a non-3-connected planar graph regularly covering a graph $H$ with a regular covering
projection $p : G \to H$. The proof is done by induction according to the size of $G$~\cite[Section
3]{negami}. Negami deals with the case that $G$ contains a 2-cut $\{u,v\}$ while stating that 1-cuts
can be solved analogously. Since $G$ is a simple graph, an inclusion minimal subgraph separated by a
2-cut is considered, giving a proper atom $A$ in $G$.  The property that $A \cap \pi(A) = \bo A \cap
\bo \pi(A)$ for every $\pi \in \Aut(G)$ of Lemma~\ref{lem:atom_automorphisms} is derived~\cite[p.
162]{negami}. 

Since non-trivial 2-cuts are not used in~\cite[p. 162]{negami}, similarly as in
Lemma~\ref{lem:proper_atoms}, $A^+$ is either 3-connected, or $K_3$. While the latter possibility is
ignored in the proof, it can be easily solved as well. We note that when trivial 2-cuts are
admitted, it is possible that $\int A \cap \pi(A) \ne \emptyset$. It happens when $G$ is a cycle,
$A$ is a subpath of length two and $\pi$ rotates vertices of $G$ by one. This is not a problem in
Negami's proof since a more restrictive version of regular graph covering is used in which the
semiregular subgroup $\Gamma$ of $\Aut(G)$ fixes no edges and further, to avoid loops in $H$, no two
adjacent vertices belong to one orbit; so the aforementioned $\pi \notin \Gamma$. Since we have no
such assumption for action of $\Gamma$, 2-cuts are required to be non-trivial.

Similarly as in Lemmas~\ref{lem:atom_covering_cases} and~\ref{lem:planar_proper_half_quotients}, it
is proved that $p |_A$ is either an edge-, a loop-, or a half-projection~\cite[p. 163, 166]{negami}
and that in all these cases $p(A)$ is planar~\cite[p. 166]{negami}. (Since no edges are fixed, only
the rotational half-quotient $p(A)$ of Lemma~\ref{lem:planar_proper_half_quotients} may appear.)

The proof is concluded by modifying $G$ and $H$ into $G'$ and $H'$ as follows. When $p |_A$ is an
edge-projection, the interiors of $p(A)$ in $H$ and of its pre-images in $G$ are either replaced by
edges, or removed when such an edge is not already in $H$. When $p |_A$ is a loop-, or a
half-projection, all these interiors are removed in $H$ and $G$. By modifying the covering
projection $p : G \to H$ into an according covering projection $p' : G' \to H'$, we get that $G'$
regularly covers $H'$.  Since $G'$ is also planar and either 3-connected, or smaller than $G$, from
the induction hypothesis every regular quotient $H'$ of $G'$ has a planar or a projectively planar
embedding. By replacing the inserted edge in this embedding with the interior of $p(A)$ (or
attaching it along an edge) or attaching the interior to a vertex in this embedding, we obtain a
planar or a projectively planar embedding of $H'$.

Suppose that we would like to get information about all regular quotients of a planar graph $G$
from~\cite{negami} as in the proof of Theorem~\ref{thm:negami}. When $G$ is a 3-connected planar
graph, all the semiregular subgroups of automorphisms can be easily
identified. In particular, $\Aut(G)$ and hence any its subgroup acts
semiregularly on the set of darts. Therefore it is enough to identify the
subgroups with trivial vertex-stabilizers.

 When $G$ is a non-3-connected planar graph, we would
like to locate proper atoms $A$ repeatedly and replace them in $G$ till we reach a 3-connected
graph. The issue is that when $p |_A$ is a loop- or a half-projection, we remove the interiors in
$G'$. This completely changes the structure of all possible quotients $G' / \Gamma'$; they do not
have to correspond to quotients $G / \Gamma$. The second issue is that we do not know $H$
beforehand, so we do not know what $p |_A$ is. The last issue is that even when new edges are
introduced, when the reduction is repeated, these added edges might be exchanged with non-added
edges in $G'$, creating new possible quotients, not corresponding to quotients of $G$.

So while Negami's paper~\cite{negami} was one of the starting points of our investigation, the
theory developed in this paper solves, among other, all the aforementioned issues with the reduction
and expansion procedures working with colored edges, loops, and half-edges of different types, as
described in Section~\ref{sec:reduction_and_expansion}.  It is key that these procedures behave well
with respect to changes in automorphism groups (Proposition~\ref{prop:reduction_homomorphism}) and
regular quotients (Lemma~\ref{lem:group_expansion}, Corollary~\ref{cor:all_group_expansions},
Theorem~\ref{thm:quotient_expansion}). Also even if Negami's more restricted definition of regular
covering and regular quotients without loops and half-edges was used, we still need to use the more
general definition in the process of the reduction, as illustrated in Fig.~\ref{fig:reducing_atoms}.

\section{Concluding Remarks} \label{sec:conclusions}

We summarize the main points addressed in this paper:
\begin{packed_itemize}
\item We describe the reduction series $G = G_0,\dots,G_r$ such that $G_{i+1}$ is constructed from
$G_i$ by replacing the atoms of $G_i$ with colored edges and the primitive graph $G_r$ is either
essentially 3-connected, or essentially $K_2$, or a essentially cycle
(Lemma~\ref{lem:primitive_graphs}). We show that $\Aut(G_i)$ is an extension of $\Aut(G_{i+1})$
(Proposition~\ref{prop:group_extension}). Changes in the automorphism groups are further studied and
applied to planar graphs in~\cite{knz}.
\item For a prescribed quotient $H_r = G_r / \Gamma_r$ of the primitive graph, we describe all
possible quotient expansions $H_0 = G_0 / \Gamma_0$ which revert the reductions.
Theorem~\ref{thm:quotient_expansion} states that every quotient $H \cong G / \Gamma$ can be obtained
in this way, and different quotients $H_0$ are constructed by non-isomorphic quotients $H_r$ and
non-isomorphic choices of half-quotients in the expansions.
\item Since the quotients of 3-connected planar graphs can be understood using geometry, we are able
to give a proof of Negami's Theorem~\cite{negami} (Theorem~\ref{thm:negami}). The reason is that
$G_r$ has a geometric embedding into the unit sphere, and a semiregular group $\Gamma_r$ can be
viewed as a spherical group acting on this embedding. Then the quotient $H_r = G_r / \Gamma_r$ is
due this geometrical interpretation embedded into the sphere, or into the projective plane. In
particular, $H_r$ is planar or projective planar. By Theorem~\ref{thm:quotient_expansion} and
Lemma~\ref{lem:planar_proper_half_quotients}, the expansions create $H$ from $H_r$ while preserving
the underlying surface of $H_r$.
\item Our results have as well algorithmic implications for regular covering testing, described
in~\cite{fkkn}. In particular, this allows to construct an algorithm for testing whether an input
planar graph $G$ regularly covers an input graph $H$, running in time $\O(\bv(G)^c \cdot 2^{\be(H)/2})$.
\end{packed_itemize}

\heading{Remarks for Reduction.} When semiregular subgroups are studied, e.g., in~\cite{knz}, it is
natural not to use halvable edges and halvable atoms. Instead, we only need to distinguish symmetric
atoms (when $\pi \in \Autbo{A}$ exchanging $\bo A$ exists) and asymmetric atoms (when $\pi$ does not
exist).

We define the reduction to  replace simultaneously the interiors of all atoms of $G_i$ by edges since
it was the simplest definition. If we replace atoms one by one,
Proposition~\ref{prop:reduction_homomorphism} does not hold. It would be possible to consider modifications   of the reduction replacing in one step just
one orbit of atoms, or just one isomorphism class of atoms.  A more interesting alternative is to replace all atoms
of only one type. For instance, if $G_i$ contains a proper atom, we only replace all proper atoms.
Otherwise if $G_i$ contains a dipole, we only replace all dipoles, and otherwise we replace all
block atoms.  The first advantage of this approach is that all edges of a dipole are original edges
of $G$ or correspond to proper atoms, similarly as all pendant edges of a star block atom are
original pendant edges or correspond to non-star block atoms. The second advantage is that the
central block/articulation is preserved in the process of reductions, without the assumption of a
non-trivial semiregular automorphism of Lemma~\ref{lem:preserved_center}.  Also, we could ignore the
reduction series $G = G_0,\dots,G_r$ and work with the reduction tree instead; see~\cite{knz} for
more details describing how $\Aut(G)$ is captured by this tree.

\heading{More General Graphs.} Our structural results also work for more general graphs. We have
assumed that the input graphs $G$ and $H$ are without loops and half-edges. We can reduce loops and
half-edges in $G$ and replace them by pendant edges. Since we assume that $H$ contains no
half-edges, we set the reductions and expansions in the way that half-edges can appear in the
expansion series but no expanded quotient $H_0$ contains half-edges. This is done by having all
edges of $G_0$ as undirected edges. To admit quotients $H_0$ with half-edges, it is sufficient to
change all edges of $G_0$ to halvable edges. Also, all the results can be used when $G$ and $H$
contain colored edges, vertices, some edges oriented, etc.

\heading{Harmonic Regular Covers.} There is a generalization of regular graph covering, admitting
singular points both in vertices and in the centers of edges, for which it would be interesting to
find out whether our techniques can be modified.  Consider geometric regular covers of surfaces,
like in Fig.~\ref{fig:quotients_of_primitive_graphs} and~\ref{fig:halfquotients_of_proper_atoms}.
The orbits of the $180^\circ$ rotations are of size two, with the exception of two points lying on
the axis of the rotation. These exceptional points are called branch points. In general, a
regular covering projection is locally homeomorphic around a branch point to the complex mapping $z
\mapsto z^\ell$ for some integer $\ell \le k$, and $\ell$ is called the \emph{order} (or index) of the branch 
point. For more details about branch points, see~\cite[Section 2.VI]{vca}.

Assume that $G$ is a 3-connected planar graph embedded onto the sphere, $\Gamma \le \Aut(G)$ is a
semiregular subgroup of automorphisms of the sphere, and $p : G \to H = G / \Gamma$ is the regular
covering projection. When $H$ is a standard graph (with no half-edges), all branch points of
$p$ belong to faces of the embedding. If a branch point (of order two) is placed in the center of an
edge of $G$, this edge is projected to a half-edge in $H$. It is possible to consider covering
projections between surfaces induced by actions of \emph{harmonic} groups, where branch points can be placed in vertices of $G$ which gives
\emph{regular harmonic  covering}~\cite{harmonic_covers}. A subgroup $\Gamma\leq \Aut(G)$ is harmonic, if it is semiregular on darts of $G$.
If a branch point of order $\ell$ is placed
in a vertex $v \in \bV(G)$, then the vertex $p(v) \in \bV(H)$ has the degree equal $\deg(v)/\ell$ and
for a dart $d \in \bD(H)$ incident with $p(v)$, the fiber $p^{-1}(d)$ has exactly $\ell$ darts
incident with $v$. In particular, every orientation preserving group of automorphisms of a 3-connected planar graph is harmonic.

\heading{4-connected Reduction.} We have described the way how to reduce a graph to a 3-connected
one while preserving its essential structural information. This approach is highly efficient for
planar graphs since many problems are much simpler for 3-connected planar graphs; for instance the
considered regular graph covering problem. Suppose that we would like to push our results further, say to
toroidal or projective planar graphs. The issue is that 3-connectivity does not restrict them much.
Is it possible to apply some ``4-connected reduction'', to reduce the input graphs even further? Suppose
that one would generalize proper atoms to be inclusion minimal parts of the graph separated by a
3-cut. Would it be possible to replace them by triangles?

\heading{Acknowledgement.}
We want to thank anonymous reviewers for many great comments which significantly improved
readability of this paper.

\bibliographystyle{plain}
\bibliography{algorithmic_regular_covers_I}

\end{document}